\documentclass{article}

\RequirePackage[T1]{fontenc}
\RequirePackage[stretch=10]{microtype}
\RequirePackage{amsmath,amsthm,amssymb}
\RequirePackage{enumerate}
\RequirePackage{graphicx}
\RequirePackage{geometry}
\RequirePackage{titlesec}
\RequirePackage[colorlinks,linktoc=page,linkcolor=blue]{hyperref}
\RequirePackage[capitalise,sort]{cleveref}
\RequirePackage{bbm}
\RequirePackage{cite}

\Crefname{assumption}{Assumption}{Assumptions}
\Crefname{figure}{Figure}{Figures}

\makeatletter

\def\Headings#1#2{\def\ps@mypagestyle{\let\@mkboth\@gobbletwo%
\def\@oddhead{\hfill {\small\scshape #1} \hfill}%
\def\@oddfoot{\hfill \small\rmfamily \thepage \hfill}%
\def\@evenhead{\hfill {\small\scshape #2} \hfill}%
\def\@evenfoot{\hfill \small\rmfamily \thepage \hfill}}%
\pagestyle{mypagestyle}}

\renewcommand\footnoterule{\kern-3\p@ \hrule \@width \textwidth \kern 2\p@}

\renewcommand{\p@enumii}{\theenumi.}
\renewcommand{\p@enumiii}{\theenumi.\theenumii.}

\makeatother

\renewenvironment{abstract}
{\centerline{\bfseries Abstract}\vspace{0.7ex}%
	\bgroup\leftskip 40pt\rightskip 40pt\small\noindent\ignorespaces}%
{\par\egroup\vskip 0.7ex}

\titleformat{\section}
{\large\bfseries}
{\thesection.}
{2mm}
{}

\titleformat{\subsection}
{\normalsize\bfseries}
{\thesubsection}
{3mm}
{}

\renewenvironment{appendix}{
\titleformat{\section}
{\large\bfseries}
{Appendix \thesection.}
{2mm}
{}
\begin{oldappendix}
}{\end{oldappendix}}

\newcommand{\Pf}[1]{{\bfseries Proof of \cref{#1}.}}

\let\oldproofname=\proofname
\renewcommand{\proofname}{\bfseries\oldproofname}

\numberwithin{equation}{section}
\newtheorem{theorem}{Theorem}[section]
\newtheorem{corollary}[theorem]{Corollary}
\newtheorem{example}[theorem]{Example}
\newtheorem{lemma}[theorem]{Lemma}
\newtheorem{proposition}[theorem]{Proposition}
\theoremstyle{definition}
\newtheorem{assumption}[theorem]{Assumption}
\newtheorem{definition}[theorem]{Definition}
\newtheorem{remark}[theorem]{Remark}

\newcommand{\B}{\mathbb{B}}

\newcommand{\N}{\mathbb{N}}

\newcommand{\R}{\mathbb{R}}

\newcommand{\Pb}{\mathbb{P}}

\newcommand{\Sb}{\mathbb{S}}

\newcommand{\Ac}{\mathcal{A}}
\newcommand{\Bc}{\mathcal{B}}
\newcommand{\Cc}{\mathcal{C}}

\newcommand{\Ec}{\mathcal{E}}
\newcommand{\Fc}{\mathcal{F}}

\newcommand{\Nc}{\mathcal{N}}

\newcommand{\Pc}{\mathcal{P}}

\renewcommand{\epsilon}{\varepsilon}

\newcommand{\Ind}{\mathbbm{1}}

\newcommand{\abs}[1]{\left|#1\right|}
\newcommand{\norm}[1]{\left\|#1\right\|}

\newcommand{\Exp}[1]{\mathbb{E}\left[#1\right]}

\newcommand{\dout}{o}
\newcommand{\dWiener}{d}
\newcommand{\Einf}[1]{\Sb_{#1 d}}
\newcommand{\CAd}[2]{\Sb^{#2}(\R^{#1})}
\newcommand{\normCAd}[2]{\norm{#1}_{\Sb^{#2}}}

\newcommand{\fct}{h}
\newcommand{\gct}{H}

\newcommand{\tr}{\mathrm{Tr}}
\newcommand{\flr}[1]{\lfloor #1 \rfloor}

\newcommand{\ED}[1]{\hat{#1}}
\newcommand{\WED}[1]{\widehat{#1}}

\newcommand{\GrowthSingle}[1]{\mathcal{G}(#1)}
\newcommand{\GrowthSum}[3]{\mathcal{G}_{#3,#2}(#1)}
\newcommand{\GrowthHoelder}[3]{\mathcal{H}_{#3,#2}(#1)}

\newcommand{\SDE}{X_{\mu,\sigma}}
\newcommand{\SDEED}{X_{\ED{\mu},\ED{\sigma}}}

\newcommand{\Euler}{\Ec^N_{\mu,\sigma}}
\newcommand{\EulerPert}{\Ec^N_{\nu,\tau}}
\newcommand{\EulED}{\Ec^N_{\ED{\mu},\ED{\sigma}}}

\newcommand{\GO}{L}
\newcommand{\GT}{l}
\newcommand{\GenO}[1]{\GO_{\mu,\sigma}^{#1}}
\newcommand{\GenT}[1]{\GT_{\mu,\sigma}^{#1}}
\newcommand{\GenDiffO}[1]{L_{\ED{\mu},\ED{\sigma}}^{#1}}
\newcommand{\GenDiffT}[1]{l_{\ED{\mu},\ED{\sigma}}^{#1}}
\newcommand{\GenBO}[1]{\mathbf{\GO}_{\mu,\sigma}^{#1}}
\newcommand{\GenBT}[1]{\mathbf{\GT}_{\mu,\sigma}^{#1}}
\newcommand{\GenStarO}[1]{\GO_{\mu_{\star},\sigma_{\star}}^{#1}}
\newcommand{\GenStarT}[1]{\GT_{\mu_{\star},\sigma_{\star}}^{#1}}

\title{\Large{\bfseries{Efficient Sobolev approximation of linear parabolic PDEs \\ in high dimensions}}\vskip 1mm}

\author{Patrick Cheridito\footnote{Department of Mathematics, ETH Zurich, Switzerland} \qquad Florian Rossmannek$^*$}

\date{}

\Headings{Efficient Sobolev approximation of linear parabolic PDEs in high dimensions}{Cheridito, Rossmannek}

\begin{document}

\maketitle

\begin{abstract}
	In this paper, we study	the error in first order Sobolev norm in the approximation of solutions to linear parabolic PDEs.
We use a Monte Carlo Euler scheme obtained from combining the Feynman--Kac representation with a Euler discretization of the underlying stochastic process.
We derive approximation rates depending on the time-discretization, the number of Monte Carlo simulations, and the dimension.
In particular, we show that the Monte Carlo Euler scheme breaks the curse of dimensionality with respect to the first order Sobolev norm.
Our argument is based on new estimates on the weak error of the Euler approximation of a diffusion process together with its derivative with respect to the initial condition.
As a consequence, we obtain that neural networks are able to approximate solutions of linear parabolic PDEs in first order Sobolev norm without the curse of dimensionality if the coefficients of the PDEs admit an efficient approximation with neural networks.
\end{abstract}


\section{Introduction}

We consider the following linear parabolic partial differential equation (PDE);
\begin{equation}
\label{Intro_PDE}
\begin{split}
	\partial_t u(t,x) + \nabla u(t,x)^T \mu(t,x) + \frac{1}{2} \tr\big( \sigma(t,x) \sigma(t,x)^T \nabla^2u(t,x) \big) + g(t,x) &= 0 \quad \text{on } [0,T] \times \R^d, \\
	u(T,\cdot) &= f \quad \text{on } \R^d.
\end{split}
\end{equation}%
There is an abundance of results that study numerical approximations of solutions of PDEs of the form \eqref{Intro_PDE}.
In physics, engineering, or finance, such equations can appear in very high dimensions.
Therefore, it is important whether numerical approximations are affected by the curse of dimensionality, meaning that the complexity of a numerical scheme grows exponentially in the dimension $d$, rendering it unfeasible in such applications.
Furthermore, in many cases, a good approximation of the gradient $\nabla u$ is just as important as a good approximation of the solution $u$.
For instance, if $u$ is the value of a financial contract, $\nabla u$ describes the sensitivities, which are intimately related to the optimal hedging strategy.

Suppose the complexity of a numerical approximation $\phi_K$ of $u$ is given by some parameter $K \in \N$.
Typical error estimates are of the form $\norm{u-\phi_K} \leq C K^{-\alpha}$ for some constants $C,\alpha > 0$.
The notion `curse of dimensionality' is not used consistently in the literature and depends on the application.
Here, we call $\phi_K$ free of the curse of dimensionality if $\alpha$ does not depend on the dimension and if the constant $C$ grows at most polynomially in the dimension.

Numerical simulations have empirically verified the ability to approximate PDEs without the curse of dimensionality, particularly machine learning paradigms;
\cite{ChanWaiNamMikaelWarin2019,BeckBeckerGrohsJaaJen2021,
EHanJen2021,HurePhamWarin2020}.
There is an ever-growing number of articles backing these simulations with proofs for breaking the curse of dimensionality; the reader may consult
\cite{BerDabGrohs2020,DeRyckMishra2021,EHanJen2021,GononSchwab2021FaS,
GrohsHerrmann2021IMANum,HutzJenKruseNguyen2022JNumM,KutyPeteRasSchnei2022}
and the references therein for a comprehensive overview.
Almost all of these results measure the error $\norm{u-\phi_K}$ in an $L^p$-norm (we comment on exceptions later).
However, in the context of PDEs, it would be more natural to also study the error incurred by derivatives.
In general, a good approximation of $u$ in $L^p$ does not have to be a good approximation in Sobolev spaces.
In applications, approximations in Sobolev spaces are more desired.
Recently, training algorithms for machine learning paradigms have been adapted to incorporate `Sobolev training';
\cite{HugeSavine2020,HurePhamWarin2020,NguThienTranCong1999,
SonJangHanHwang2021,Tsay2021,VlassisSun2021,YuanZhuLiuDingZhangZhang2022}.
In this Sobolev training, an algorithm takes into account loss values of the derivative approximation of a model;
\cite{CocolaHand2020,CzarOsinJaderSwirPasc2017,KisselDiepold2020}.
Physically inspired neural networks are trained so that they satisfy a given PDE as well as possible;
\cite{DeRyckMishra2021,SiriSpili2018}.
In this approach, derivatives are trained implicitly.
In terms of approximation capacities, the theoretical justification for these algorithms is insufficient.
They are mainly based on the ability of a model class to admit suitable approximations of a target function in Sobolev space;
\cite{AbdeljawadGrohs2022,DeRyckLanthMishra2021,GuehKutyPete2020,
HonYang2022,HorStinWhite1990,RolTeg2018,SiriSpili2018}.
But these are usually not quantitative or at least suffer the curse of dimensionality (we comment on exceptions later).
In this article, we improve the theoretical justification for the Sobolev training of PDEs by proving error rates for approximating $u$ and its first derivative without the curse of dimensionality.

A common approach for approximating solutions to linear parabolic PDEs is to exploit the Feynman--Kac formula and move the approximation problem to the realm of stochastic processes.
For a simpler presentation, we here assume $g=0$ but treat non-zero $g$ later.
Omitting technical assumptions, the Feynman--Kac formula states that
\begin{equation}
\label{Intro_FeynmanKac}
	u(0,x) = \Exp{ f(X_T) },
\end{equation}%
where $u$ is the solution to PDE \eqref{Intro_PDE} and $X$ is the solution to the corresponding stochastic differential equation (SDE);
\begin{equation}
\label{Intro_SDE}
\begin{split}
	dX_t &= \mu(t,X_t) dt + \sigma(t,X_t) dW_t,
	\\
	X_0 &= x,
\end{split}
\end{equation}%
where $W$ is a standard $d$-dimensional Brownian motion.
The functions $\mu$ and $\sigma$ are called drift and diffusion coefficient, respectively.
The expected value in the Feynman--Kac formula can be tackled with a Monte Carlo simulation, which is famously known not to suffer the curse of dimensionality.
Thus, the problem reduces to efficiently approximating the stochastic process $X$.
Many numerical schemes have been studied for this purpose, the most prominent one being the Euler scheme and variations of it.
Fix a deterministic time grid $0 = t_0 < t_1 < \dots < t_N = T$.
The explicit Euler scheme $\Ec$ in its standard form is given by
\begin{equation}
\label{Intro_Euler}
\begin{split}
	\Ec_{t_{n+1}} &= \Ec_{t_n} + \mu(t_n,\Ec_{t_n}) (t_{n+1}-t_n) + \sigma(t_n,\Ec_{t_n}) (W_{t_{n+1}}-W_{t_n}),
	\\
	\Ec_0 &= x.
\end{split}
\end{equation}%
We now state an informal version of our main result.

\begin{theorem}[informal]
\label{Intro_theorem}
	Let $\mu$, $\sigma$, $f$, $g$, and $u$ adhere to certain regularity assumptions, where $u$ is a solution to PDE \eqref{Intro_PDE}.
Let $\Ec_1(x),\dots,\Ec_M(x)$ be $M \in \N$ independent copies of the Euler scheme \eqref{Intro_Euler} for SDE \eqref{Intro_SDE} based on a uniform time grid with $N \in \N$ steps.
Consider $\phi \colon \Omega \times \R^d \rightarrow \R^{\dout}$ given by
\begin{equation}
\label{Intro_theorem_approximator}
\begin{split}
	\phi(x) &= \frac{1}{M} \sum_{m=1}^M \left[ f\left(\Ec_m(x)_T\right) + \frac{T}{N} \sum_{n=0}^{N-1} g\left(t_n,\Ec_m(x)_{t_n}\right) \right].
\end{split}
\end{equation}%
In the following, $c \geq 0$ denotes a constant that depends only on $\mu$ and $\sigma$, and $c_{\alpha} \geq 1$ denotes a constant that depends on $\mu$, $\sigma$, $f$, $g$, $u$, $T$, and $\alpha \in (0,1]$.
Now, the probabilistic numerical scheme $\phi$ approximates the PDE solution $U = u(0,\cdot)$ at time 0 in the following sense.
Let $\delta,\epsilon \in (0,1)$ and let $K \subseteq \R^d$ be compact.
If $N \geq 16cT$, then
\begin{equation*}
\begin{split}
	\Pb\left( \norm{U - \phi}_{L^2(K)} + \norm{\nabla U - \nabla\phi}_{L^2(K)} \leq \epsilon \right)
	&\geq 1 - \delta
\end{split}
\end{equation*}%
provided $M$ and $N$ are so large that
\begin{equation*}
	\frac{\epsilon \sqrt{\delta}}{M^{-1/2} + \sqrt{\delta} N^{-\alpha}}
	\geq c_{\alpha} d^2 \norm{1 + \norm{x}^{10}}_{L^2(K)}.
\end{equation*}%
\end{theorem}

The regularity assumptions required and the constants $c$ and $c_{\alpha}$ will be made explicit in our main result, \cref{Sobolev_MCES_error}.
We remark that these constants can depend implicitly on the dimension since the functions $\mu$, $\sigma$, $f$, $g$, and $u$ do.
The specification in \cref{Sobolev_MCES_error} unfolds this implicit dependence, linking it to the spatial growth of the functions.
In particular, it reveals in which case the curse of dimensionality is broken to the full extent.
Additionally, we will argue why we cannot expect to alleviate the implicit dependence.
Let us mention that starting the Euler scheme at time $t_0=0$ is convenient for the sake of notation.
But if we let it start at an arbitrary initial time $t \in [0,T]$ and define the Euler scheme on the uniform time grid in $[t,T]$, then we obtain a time-dependent $\phi(t,x)$ that approximates $u$ in the norm $\norm{u - \phi}_{L^2([0,T] \times K)} + \norm{\nabla u - \nabla\phi}_{L^2([0,T] \times K)}$ at an analogous rate.

There is a vast amount of results on convergence properties and error rates of the Euler scheme.
We summarize the findings most relevant to us here but refer to the exposition in
\cite{BossyJabirMartinez2021,HutzJenKloe2011}
and the references
\cite{NgoTaguchi2016MComp,BerBossyDiop2008ESAIM,
BossyJabirMartinez2021,MilTret2005,HutzJen2020,HuLiMao2018}
for a more detailed account of the existing literature.
Typically, to obtain not only mere convergence but also a rate of convergence, assumptions imposed on the drift and diffusion coefficients involve global Lipschitz continuity in the space variable and H{\"o}lder continuity in the time variable.
To relax the global Lipschitz continuity for both coefficients simultaneously, \cite{YuanMao2008SAA} instead imposed on them to grow at most linearly and be locally Lipschitz continuous with local Lipschitz constants on balls growing at most logarithmically in the diameter.
The first of these assumptions is sharp in the following sense:
\cite{HutzJenKloe2011,HutzJenKloe2013} showed that the Euler scheme diverges as soon as one of the coefficient functions grows more than linearly.
In this article, we narrow the gap to the negative results in
\cite{HutzJenKloe2011,HutzJenKloe2013}
by proving a convergence rate for the Euler scheme for at most linearly growing coefficients with a much weaker growth assumption on the local Lipschitz constants than in \cite{YuanMao2008SAA}.
Our result, \cref{weak_error_rate_thrm}, is based on a weak error estimate, which is the reason for the regularity requirements in \cref{Intro_theorem}.

Let us illustrate in more detail our method of attack to obtain Sobolev approximations of PDE \eqref{Intro_PDE} based on the Feynman--Kac formula.
A key observation is the well-known differentiability of the stochastic process $X = X(x)$ and its Euler scheme $\Ec = \Ec(x)$ with respect to the initial point $x$.
From the Feynman--Kac formula \eqref{Intro_FeynmanKac}, it is to be expected that
\begin{equation}
\label{Intro_FeynmanKac_Diff}
	\nabla u(0,x) = \Exp{ \nabla f(X_T)^T DX_T } \approx \Exp{ \nabla f(\Ec_T)^T D\Ec_T }.
\end{equation}%
The derivative $D\Ec$ of the Euler scheme $\Ec$ is given by
\begin{equation*}
	D\Ec_{t_{n+1}} = D\Ec_{t_n} + \nabla \mu(t_n,\Ec_{t_n})^T D\Ec_{t_n} (t_{n+1}-t_n) + \nabla \sigma(t_n,\Ec_{t_n})^T D\Ec_{t_n} (W_{t_{n+1}}-W_{t_n}).
\end{equation*}%
We cannot apply known error estimates for Euler schemes to $D\Ec$ because $D\Ec$ is not the Euler scheme of $DX$.
Nonetheless, we can regard $D\Ec$ as a Euler scheme by considering it simultaneously with $\Ec$.
Indeed, the Euler scheme of the combined process $(X,DX)$ is exactly $(\Ec,D\Ec)$.
The combined process $(X,DX)$ solves the SDE $dZ_t = \ED{\mu}(t,Z_t) dt + \ED{\sigma}(t,Z_t) dW_t$ with the new drift coefficient $\ED{\mu}(t,x,y) = (\mu(t,x),\nabla\mu(t,x)^Ty)$ and diffusion coefficient $\ED{\sigma}(t,x,y) = (\sigma(t,x),\nabla\sigma(t,x)^Ty)$.
Now that we can view $D\Ec$ as a Euler scheme, we can tackle the error estimate in \eqref{Intro_FeynmanKac_Diff}.
The existing tools surveyed above are not applicable to the coefficients $\ED{\mu}$ and $\ED{\sigma}$, but we will be able to apply the new convergence rate result from this article.

Let us mention a special case, in which Sobolev approximation of PDEs free of the curse of dimensionality has been achieved.
This has been done for the heat equation in \cite{BerDabGrohs2020} and, more generally, for the PDE \eqref{Intro_PDE} under the assumption that $\mu$ and $\sigma$ are affine functions of $x$ in \cite{DeRyckMishra2021}.
The key to this result is that the stochastic process $X = X(x)$ itself becomes an affine function of $x$ and the approximation error can be estimated with a simple Monte Carlo argument.
More precisely, there exist random variables $A$ and $B$ independent of $x$ such that $X_T(x) = Ax+B$.
Then, there exist matrices $A_m$ and vectors $B_m$ such that \eqref{Intro_FeynmanKac_Diff} simplifies to
\begin{equation*}
	\nabla u(0,x) = \Exp{ \nabla f(Ax+B)^T A } \approx \frac{1}{M} \sum_{m=1}^M \nabla f(A_mx+B_m)^T A_m.
\end{equation*}%
In this article, we extend the findings of
\cite{BerDabGrohs2020,DeRyckMishra2021}
to non-affine coefficients.
Lastly, we present an application to the approximation theory of neural networks.
Networks have also been employed in \cite{BerDabGrohs2020,DeRyckMishra2021}, among other works mentioned above.
This is due to the compositional structure of networks, which makes them suitable to implement the Monte Carlo Euler scheme from \cref{Intro_theorem}.
If the coefficients $\mu$ and $\sigma$ and the functions $f$ and $g$ can be represented by networks, then the function $\phi(\omega,\cdot)$ given in \eqref{Intro_theorem_approximator} can also be represented by a network element-wise for each $\omega \in \Omega$.
This gives rise to a network Sobolev approximation of the solution to PDE \eqref{Intro_PDE}.
If the coefficients cannot be represented exactly, then we approximate them with networks and build an analogous function $\phi$, where $\mu$, $\sigma$, $f$, and $g$ are replaced by their network approximations.
This yields a perturbed scheme, for which we derive corresponding error estimates.
To obtain a network approximation that breaks the curse of dimensionality, we require the approximations of the coefficients to be efficient.
This is done in \cref{NN_error}.
As such, we provide the first proof that networks can approximate solutions to PDEs with non-affine coefficients in a Sobolev sense without the curse of dimensionality.

Our exposition is structured as follows.
We prove the new weak error rate of the Euler scheme for at most linearly growing coefficients in \cref{section_weak_error}.
Thereafter, we add a global Lipschitz assumption on the coefficients to obtain a rate of convergence in the first order Sobolev norm.
The perturbed scheme is analyzed in \cref{section_perturbed}, first in general and then for networks.
In the appendix, we develop the technical `polynomial growth calculus', which proves useful for keeping track of implicit dependencies of constants on the dimension.


\subsection*{Notation}

Euclidean spaces $\R^{d \times \dWiener}$ are endowed with the Frobenius norm $\norm{A} = \sqrt{\tr(A^T A)}$.
Throughout, we fix a time horizon $T \geq 1$, measurable drift and diffusion coefficients $\mu \colon [0,T] \times \R^d \rightarrow \R^d$ and $\sigma \colon [0,T] \times \R^d \rightarrow \R^{d \times \dWiener}$, the inhomogeneity of the PDE $g \colon [0,T] \times \R^d \rightarrow \R^{\dout}$, and the terminal value of the PDE $f \colon \R^d \rightarrow \R^{\dout}$.
We let $\GenT{t}$ be the operator acting on twice differentiable functions $\fct \colon \R^d \rightarrow \R^{\dout}$ by\footnote{We understand $\tr\big( \sigma(y,t) \sigma(y,t)^T \nabla^2\fct(x) \big)$ as the $\dout$-dimensional vector with elements $\tr\big(\sigma(y,t) \sigma(y,t)^T \nabla^2\fct_i(x)\big)$, $i = 1,\dots,\dout$, for $\fct(x) = (\fct_1(x),\dots,\fct_{\dout}(x))$.}
\begin{equation*}
	\GenT{t} \fct \colon \R^d \times \R^d \rightarrow \R^{\dout}, \quad \GenT{t} \fct(x,y) = \nabla \fct(x)^T \mu(t,y) + \frac{1}{2} \tr\big( \sigma(t,y) \sigma(t,y)^T \nabla^2\fct(x) \big).
\end{equation*}%
For a time-dependent function $\fct \colon [0,T] \times \R^d \rightarrow \R^{\dout}$, we write $\GenT{t}\fct(s,x,y) = \GenT{t}[\fct(s,\cdot)](x,y)$.
Then, with the operator $\GenO{t}$ given by
\begin{equation*}
	\GenO{t} \fct \colon \R^d \rightarrow \R^{\dout}, \qquad \GenO{t} \fct(x) = \GenT{t} \fct(x,x),
\end{equation*}%
we can write PDE \eqref{Intro_PDE} succinctly as $\partial_t u + \GenO{}u + g = 0$.
We call $\GenO{}$ the generator.

For the stochastic part, we work on a complete probability space $(\Omega,\Fc_0,\Pb)$ with a filtration $\Fc = (\Fc_t)_{t \geq 0}$ satisfying the usual conditions and a standard $\R^{\dWiener}$-valued Wiener process $(W_t)_{t \geq 0}$ for $\Fc$.
The norm on $L^p(\Omega,\R^d)$, $p \in [2,\infty)$, will be denoted by $\norm{X}_p = \Exp{ \norm{X}^p }^{1/p}$.
We let $(\CAd{d}{p},\normCAd{\cdot}{p})$ be the Banach space $\CAd{d}{p} = \left\{ X \in L^p\big(\Omega,C([0,T],\R^d)\big) \colon X \text{ is } \Fc\text{-adapted} \right\}$ of all $\Fc$-adapted $\Pb$-a.s.\ continuous $\R^d$-valued stochastic processes on $[0,T]$ with norm
\begin{equation*}
	\normCAd{X}{p} = \Exp{ \sup_{t \in [0,T]} \norm{X_t}^p }^{1/p}.
\end{equation*}%
Let $\Einf{} \subseteq \CAd{d}{2}$ be the subset $\Einf{} = \bigcap_{p \geq 2} \CAd{d}{p}$.


\section{Weak error rate for the Monte Carlo Euler scheme}
\label{section_weak_error}

\subsection{Measuring polynomial growth}

As discussed in the introduction, breaking the curse of dimensionality in an approximation task involves polynomial control on the constant $C$ in the error estimate with respect to the dimension.
Any reasonable numerical scheme for approximating solutions of PDEs of the form \eqref{Intro_PDE} and with it its error estimate will in one way or another depend on the drift and diffusion coefficients.
The spatial domain of $\mu$ and $\sigma$ is $\R^d$, which unveils an implicit dependence on the dimension.
Thus, we need to quantify how the approximation error depends on $\mu$ and $\sigma$ to guarantee that this implicit dependence does not incur the curse of dimensionality.
We approach this with measuring the polynomial growth of a function.
Since we study Sobolev approximations and will rely on a weak error estimate for the Euler scheme, we consider the growth of derivatives.

\begin{definition}
\label{def_PG}
	For a $\N_0 \ni k$-times differentiable function $\fct \colon \R^d \rightarrow \R^{\dout}$ and $r \in [0,\infty)$, we introduce the growth measure
\begin{equation*}
	\GrowthSum{\fct}{r}{k} = \sup_{x \in \R^d} \sum_{i=0}^k \frac{\norm{\nabla^i \fct(x)}}{(1+\norm{x})^{r+k-i}} \in [0,\infty].
\end{equation*}%
If $\fct \colon [0,T] \times \R^d \rightarrow \R^{\dout}$ depends on time, then we let $\GrowthSum{\fct}{r}{k} = \sup_{t \in [0,T]} \GrowthSum{\fct(t,\cdot)}{r}{k}$.
Since $\GrowthSum{\cdot}{1}{0}$ will appear frequently, we shorten this to
\begin{equation*}
	\GrowthSingle{\fct}
	= \GrowthSum{\fct}{1}{0}
	= \sup_{x \in \R^d} \frac{\norm{\fct(x)}}{1+\norm{x}}.
\end{equation*}%
\end{definition}

\begin{example}
\label{ex_PG}
	Let $\fct \colon \R^d \rightarrow \R^{\dout}$ be Lipschitz continuous with Lipschitz constant $L \geq 0$.
Then, for all $r \in [1,\infty)$,
\begin{equation*}
	\GrowthSum{\fct}{r}{0}
	\leq \max\left\{ \norm{\fct(0)} , \frac{L}{r} \right\}.
\end{equation*}%
If, in addition, $\fct$ is differentiable, then, for all $r \in [0,\infty)$,
\begin{equation*}
	\GrowthSum{\fct}{r}{1}
	\leq \sqrt{\min\{ d , \dout \}} L + \max\left\{ \norm{\fct(0)} , \frac{L}{r+1} \right\}.
\end{equation*}%
\end{example}

We remark that if $\GrowthSum{\nabla^k \fct}{r}{0} < \infty$, then $\GrowthSum{\fct}{r}{k} < \infty$; see \cref{PGcalc_diff} in the appendix.
The converse is obvious.
The exact definition of $\GrowthSum{\fct}{r}{k}$ with the power $r+k-i$ in the denominator is for convenience, because it allows for concise inequalities such as the following.

\begin{lemma}
\label{PGcalc_int}
	Let $\fct \colon \R^d \rightarrow \R^{\dout}$ be $\N_0 \ni k$-times differentiable and $r \in [0,\infty)$.
Then, for any $0 \leq l,m \leq k$ with $l+m \leq k$, we have
\begin{equation*}
	\GrowthSum{\nabla^l \fct}{r+k-l-m}{m}
	\leq \GrowthSum{\fct}{r}{k}.
\end{equation*}%
\end{lemma}

\begin{proof}
	This is immediate from a change of summation index in the definition of $\GrowthSum{\nabla^l \fct}{r+k-l-m}{m}$.
\end{proof}

We will require the drift and diffusion coefficients to be H{\"o}lder continuous in time, but allow the H{\"o}lder constant to grow in space.

\begin{definition}
	For a function $\fct \colon [0,T] \times \R^d \rightarrow \R^{\dout}$ and $\alpha \in (0,1]$, $r \in [0,\infty)$, we introduce the H{\"o}lder growth measure
\begin{equation*}
		\GrowthHoelder{\fct}{r}{\alpha} = \sup_{\substack{s,t \in [0,T] \\ 0 < |t-s| \leq 1}} \sup_{x \in \R^d} \frac{\norm{\fct(t,x)-\fct(s,x)}}{(1+\norm{x})^r} |t-s|^{-\alpha} \in [0,\infty].
\end{equation*}%
\end{definition}

The reason for introducing these growth measures is to develop a `polynomial growth calculus', with which we can conveniently estimate the growth of a function by that of related functions without knowing anything else about those other functions.
An example of this was given in \cref{PGcalc_int}.
Another example, which we will use later, is the following estimate for the growth of the generator applied to a function twice.
To state this example, it is convenient to introduce the operators $\GenBO{t}$ and $\GenBT{t}$ acting on twice differentiable functions $\fct \colon \R^d \times \R^d \rightarrow \R^{\dout}$ by
\begin{equation*}
\begin{alignedat}{2}
	\GenBT{t}\fct \colon \R^d \times \R^d \times \R^d \rightarrow \R^{\dout},& \qquad &&\GenBT{t}\fct(x,y,z) = \GenT{t}[\fct(\cdot,z)](x,y),
	\\
	\GenBO{t}\fct \colon \R^d \times \R^d \rightarrow \R^{\dout},& \qquad &&\GenBO{t}\fct(x,z) = \GenBT{t}\fct(x,x,z) = \GenO{t}[\fct(\cdot,z)](x).
\end{alignedat}%
\end{equation*}%

\begin{lemma}
\label{PGcalc_double_generator}
	Assume $\mu$ and $\sigma$ are twice differentiable in the space variable and let $r_1,r_2,r_3 \in [0,\infty)$.
\begin{enumerate}[\rm (i)]\itemsep = 0em
\item
For any twice differentiable function $\fct \colon \R^{2d} \rightarrow \R^{\dout}$ and any $t \in [0,T]$,
\begin{equation}
\label{PGcalc_double_generator_estimate_1}
	\GrowthSum{\GenBO{t}\fct}{r}{0}
	\leq \GrowthSum{\GenBT{t}\fct}{r}{0}
	\leq \max\left\{ \GrowthSum{\mu}{r_1}{0} , \frac{1}{2} \GrowthSum{\sigma}{r_2}{0}^2 \right\} \GrowthSum{\nabla \fct}{r_3}{1}
\end{equation}%
with $r = \max\{r_1 + 1, 2r_2\} + r_3$.

\item
For any four times differentiable function $\fct \colon \R^d \rightarrow \R^{\dout}$ and any $t_1,t_2 \in [0,T]$,
\begin{equation*}
\begin{split}
	\GrowthSum{\GenBT{t_2}\GenT{t_1}\fct}{r}{0}
	&\leq 2^{r+r_3+14} d \max\left\{ \GrowthSum{\mu}{r_1}{2}^2 , \frac{1}{4} \GrowthSum{\sigma}{r_2}{2}^4 \right\} \GrowthSum{\fct}{r_3}{4}
\end{split}
\end{equation*}%
with $r = 2 \max\{r_1, 2r_2 + 1\} + r_3 + 8$, and the same bound holds for $\GrowthSum{\GenT{t_2}\GenO{t_1}\fct}{r}{0}$.
\end{enumerate}
\end{lemma}

Estimate \eqref{PGcalc_double_generator_estimate_1} also holds for $\GenO{t}\fct$ and $\GenT{t}\fct$ in place of $\GenBO{t}\fct$ and $\GenBT{t}\fct$ for a function $\fct \colon \R^d \rightarrow \R^{\dout}$ with domain $\R^d$, since we can apply the lemma with $(w,z) \mapsto \fct(w)$.
The proof of \cref{PGcalc_double_generator} is deferred to the appendix, in which we develop a more elaborate `polynomial growth calculus'.

There is one particular cumbersome constant that will appear in several of the results to come.
To avoid repeating it throughout, we introduce it here as
\begin{equation}
\label{cumbersome_constant}
\begin{split}
	\Cc_{r,\alpha}(\mu,\sigma,g,u)
	&= \max\big\{ \GrowthSum{\partial_t g}{r}{0} , \max\left\{ \GrowthSum{\mu}{r}{2} , \GrowthSum{\sigma}{r}{2}^2 \right\} \GrowthSum{g}{r}{2} ,
	\\
	& \qquad\qquad \max\left\{ \GrowthSum{\mu}{r}{2}^2 , \GrowthSum{\sigma}{r}{2}^4 , \GrowthHoelder{\mu}{r}{\alpha} , \GrowthSum{\sigma}{r}{2} \GrowthHoelder{\sigma}{r}{\alpha} \right\} \GrowthSum{u}{r}{4} \big\}.
\end{split}
\end{equation}%


\subsection{Weak error rate for the Euler scheme}

The Euler scheme is usually regarded as a temporal discretization of an SDE.
However, for a weak error analysis, one can introduce the Euler scheme directly without considering the original SDE.
Fix a natural number $N \in \N$ and abbreviate $t_n = nT/N$.
Given $t \in [0,T]$, we write $\flr{t} = t_n$ for $t \in [t_n,t_{n+1})$, $0 \leq n \leq N$.
Then, for any $x \in \R^d$, we let $\Ec^N(x)$ be the Euler scheme
\begin{equation*}
\begin{split}
	\Ec^N(x)_t &= x + \int_0^t \mu(\flr{s},\Ec^N(x)_{\flr{s}}) ds + \int_0^t \sigma(\flr{s},\Ec^N(x)_{\flr{s}}) dW_s, \qquad 0 \leq t \leq T,
	\\
	\Ec^N(x)_0 &= x.
\end{split}
\end{equation*}%
Here, we take $t_0 = 0$ as the initial time to keep the notation simple, but all arguments go through for an arbitrary initial time $t \in [0,T]$.
It is not difficult to show that the Euler scheme $\Ec^N(x)$ is an element of $\Einf{}$ if $\GrowthSum{\mu}{r}{0},\GrowthSum{\sigma}{r}{0} < \infty$ for some $r \geq 0$.
In the lemma below, we derive a bound on $\norm{\Ec^N(x)_t}_p$ in the case of at most linearly growing $\mu$ and $\sigma$, that is $r=1$.
The qualitative version of this lemma is well-known but we need a quantitative version with explicit constants.

\begin{lemma}
\label{Euler_scheme_well_defined}
	Let $p \in [2,\infty)$ and $c = \max\{\GrowthSingle{\mu},\GrowthSingle{\sigma}^2\}$.
If $N \geq 16cT$, then there exists a constant $\kappa_p > 0$ that depends only on $p$ such that for all $x \in \R^d$
\begin{equation*}
	\sup_{t \in [0,T]} \norm{\Ec^N(x)_t}_p^p
	\leq e^{\kappa_p c T} (\kappa_p c T+\norm{x}^p).
\end{equation*}%
\end{lemma}

\begin{proof}
	Throughout, $\kappa \geq 1$ denotes a constant that depends only on $p$ but may change from line to line.
Abbreviate $\Ec = \Ec^N(x)$.
By It{\^o}'s formula applied to the function $\fct(x) = \norm{x}^p$, we have
\begin{equation}
\label{Pf_Eswd_Ito}
	\Exp{ \norm{\Ec_t}^p } = \norm{x}^p + \Exp{ \int_0^t \GenT{\flr{s}} \fct(\Ec_s,\Ec_{\flr{s}}) ds }.
\end{equation}%
For any $x,y \in \R^d$,
\begin{equation*}
\begin{split}
	|\GenT{\flr{s}} \fct(x,y)|
	&\leq p \norm{x}^{p-1} \norm{\mu(\flr{s},y)} + p (p-1) \norm{x}^{p-2} \norm{\sigma(\flr{s},y)}^2
	\\
	&\leq c p \norm{x}^{p-1} ( 1 + \norm{y} ) + c p (p-1) \norm{x}^{p-2} ( 1 + \norm{y} )^2
	\\
	&\leq \kappa c \left( 1 + \norm{x}^p + \norm{x}^{p-1} \norm{y} + \norm{x}^{p-2} \norm{y}^2 \right).
\end{split}
\end{equation*}%
Thus, by H{\"o}lder's inequality,
\begin{equation}
\label{Pf_Eswd_Ito_integrand}
\begin{split}
	\Exp{ | \GenT{\flr{s}} \fct(\Ec_s,\Ec_{\flr{s}}) | }
	&\leq \kappa c \left( 1 + \Exp{ \norm{\Ec_s}^p } + \Exp{ \norm{\Ec_s}^{p-1} \norm{\Ec_{\flr{s}}} } + \Exp{ \norm{\Ec_s}^{p-2} \norm{\Ec_{\flr{s}}}^2 } \right)
	\\
	&\leq \kappa c \left( 1 + \Exp{ \norm{\Ec_s}^p } + \Exp{ \norm{\Ec_s}^p }^{\frac{p-1}{p}} \Exp{ \norm{\Ec_{\flr{s}}}^p }^{\frac{1}{p}} + \Exp{ \norm{\Ec_s}^p }^{\frac{p-2}{p}} \Exp{ \norm{\Ec_{\flr{s}}}^p }^{\frac{2}{p}} \right).
\end{split}
\end{equation}%
For each $n$ and $s \in [t_n,t_{n+1}]$, we can apply Doob's maximal inequality to the martingale $\Ec_r - \mu(t_n,\Ec_{t_n})(r-t_n)$ indexed by $[t_n,s]$ to find
\begin{equation*}
\begin{split}
	\Exp{ \sup_{t_n \leq r \leq s} \norm{\Ec_{r}}^p }^{\frac{1}{p}}
	&\leq \frac{p}{p-1} \Exp{ \norm{ \Ec_s - \mu(t_n,\Ec_{t_n}) (s-t_n) }^p }^{\frac{1}{p}}
	\\
	&\leq 4 \Exp{ \norm{ \Ec_s }^p }^{\frac{1}{p}} + 4 \Exp{ \norm{ \mu(t_n,\Ec_{t_n}) (s-t_n) }^p }^{\frac{1}{p}}
	\\
	&\leq 4 \Exp{ \norm{ \Ec_s }^p }^{\frac{1}{p}} + 8 c T N^{-1} \left( 1 + \Exp{ \norm{ \Ec_{t_n} }^p }^{\frac{1}{p}} \right).
\end{split}
\end{equation*}%
Hence, if $N \geq 16cT$, then
\begin{equation*}
	\Exp{ \norm{\Ec_{\flr{s}}}^p }^{\frac{1}{p}}
	\leq \Exp{ \sup_{t_n \leq r \leq s} \norm{\Ec_{r}}^p }^{\frac{1}{p}}
	\leq 1 + 8 \Exp{ \norm{ \Ec_s }^p }^{\frac{1}{p}}.
\end{equation*}%
The estimate in \eqref{Pf_Eswd_Ito_integrand} simplifies to $\Exp{ | \GenT{\flr{s}} \fct(\Ec_s,\Ec_{\flr{s}}) | } \leq \kappa c \left( 1 + \Exp{ \norm{\Ec_s}^p } \right)$.
We plug this into \eqref{Pf_Eswd_Ito};
\begin{equation*}
	\Exp{ \norm{\Ec_t}^p }
	\leq \norm{x}^p + \kappa c t + \kappa c \int_0^t \Exp{ \norm{\Ec_s}^p } ds.
\end{equation*}%
Finally, by Gr{\"o}nwall's inequality, $\Exp{ \norm{\Ec_t}^p } \leq \left( \norm{x}^p + \kappa c t \right) e^{\kappa c t}$.
\end{proof}

We have laid the groundwork for our weak error analysis of the Euler scheme in the approximation of PDE solutions.
The following smoothness assumptions are needed to apply It{\^o}'s lemma in the proof of \cref{weak_error_rate_thrm}.

\begin{assumption}
\label{assumption_weak_error}
	Suppose $\mu$, $\sigma$, and $g$ are twice continuously differentiable in the space variable, $g$ is also continuously differentiable in the time variable, and $f$ is four times continuously differentiable.
Suppose $u \in C^{1,4}([0,T] \times \R^d,\R^{\dout})$ solves the PDE
\begin{equation*}
\begin{split}
	\partial_t u + \GenO{}u + g &= 0 \quad \text{on } [0,T] \times \R^d , \\
	u(T,\cdot) &= f \quad \text{on } \R^d.
\end{split}
\end{equation*}%
\end{assumption}

We remark that growth assumptions on $\mu$ and $\sigma$ that guarantee everything to be well-posed are implicit.
For if such growth assumptions were not met, then the bound in the theorem below and in the other results to come would have infinity on the right hand side and the statement of the theorem would become void.

\begin{theorem}
\label{weak_error_rate_thrm}
	Let \cref{assumption_weak_error} hold and consider $v \colon \R^d \rightarrow \R^{\dout}$ given by
\begin{equation*}
	v(x) = \Exp{ f\left(\Ec^N(x)_T\right) + \frac{T}{N} \sum_{n=0}^{N-1} g\left(t_n,\Ec^N(x)_{t_n}\right) }.
\end{equation*}%
Let $\alpha \in (0,1]$, $r \in [0,\infty)$, and denote $c = \max\{\GrowthSingle{\mu},\GrowthSingle{\sigma}^2\}$ and $C_{r,\alpha} = \Cc_{r,\alpha}(\mu,\sigma,g,u)$; see \eqref{cumbersome_constant}.
If $N \geq 16 c T$, then there exists a constant $\kappa_r > 0$ depending only on $r$ such that for all $x \in \R^d$
\begin{equation*}
	\norm{ u(0,x) - v(x) }
	\leq d C_{r,\alpha} \kappa_r (1+c) T^3 e^{ \kappa_r c T} \left( 1 + \norm{x}^{5r+10} \right) \frac{1}{N^{\alpha}}.
\end{equation*}%
\end{theorem}

\begin{proof}
	Abbreviate $\Ec = \Ec^N(x)$.
Since $u(T,\cdot) = f$, we have
\begin{equation*}
	v(x) - u(0,x) = \Exp{u(T,\Ec_T)} - \Exp{u(0,\Ec_0)} + \frac{T}{N} \sum_{n=0}^{N-1} \Exp{g(t_n,\Ec_{t_n})}.
\end{equation*}%
By It{\^o}'s formula and the PDE solved by $u$,
\begin{equation*}
\begin{split}
	\Exp{u(T,\Ec_T)} - \Exp{u(0,\Ec_0)}
	&= \Exp{ \int_0^T - g(s,\Ec_s) - \GenO{s}u(s,\Ec_s) + \GenT{\flr{s}} u(s,\Ec_s,\Ec_{\flr{s}}) ds }.
\end{split}
\end{equation*}%
From now on, we write
\begin{equation*}
\begin{split}
	\eta_1 &= \Exp{ \int_0^T - \GenO{s}u(s,\Ec_s) + \GenT{\flr{s}} u(s,\Ec_s,\Ec_{\flr{s}}) ds },
	\\
	\eta_2 &= \frac{T}{N} \sum_{n=0}^{N-1} \Exp{g(t_n,\Ec_{t_n})} - \Exp{ \int_0^T g(s,\Ec_s) ds }
\end{split}
\end{equation*}%
so that $v(x) - u(0,x) = \eta_1 + \eta_2$.
We first bound $\eta_1$, which we split into the three terms
\begin{equation*}
\begin{split}
	\eta_1
	&= \Exp{ \int_0^T \GenO{s}u(s,\Ec_{\flr{s}}) - \GenO{s}u(s,\Ec_s) ds }
	+ \Exp{ \int_0^T \GenO{\flr{s}} u(s,\Ec_{\flr{s}}) - \GenO{s}u(s,\Ec_{\flr{s}}) ds }
	\\
	&\quad+ \Exp{ \int_0^T \GenT{\flr{s}} u(s,\Ec_s,\Ec_{\flr{s}}) - \GenT{\flr{s}} u(s,\Ec_{\flr{s}},\Ec_{\flr{s}}) ds } =: \eta_{1,1} + \eta_{1,2} + \eta_{1,3}.
\end{split}
\end{equation*}%
We can apply It{\^o}'s formula to $\GenO{s}u(s,\cdot)$ for each $s \in [0,T]$ to rewrite $\eta_{1,1}$ as
\begin{equation*}
\begin{split}
	\eta_{1,1}
	&= - \int_0^T \Exp{ \GenO{s}u(s,\Ec_s) - \GenO{s}u(s,\Ec_{\flr{s}}) } ds
	= - \int_0^T \Exp{ \int_{\flr{s}}^s \GenT{\flr{\tau}} \GenO{s}u(s,\Ec_{\tau},\Ec_{\flr{\tau}}) d\tau } ds.
\end{split}
\end{equation*}%
In the estimates to come, $\kappa \geq 1$ denotes a constant that depends only on $r$ but may change from line to line.
By \cref{PGcalc_double_generator},
\begin{equation*}
\begin{split}
	\Exp{ \norm{ \GenT{\flr{\tau}} \GenO{s}u(s,\Ec_{\tau},\Ec_{\flr{\tau}}) } }
	&\leq \GrowthSum{\GenT{\flr{\tau}} \GenO{s}u(s,\cdot)}{5r+10}{0} \Exp{ \left( 1 + \norm{(\Ec_{\tau},\Ec_{\flr{\tau}})} \right)^{5r+10} }
	\\
	&\leq \kappa d C \left( 1 + \sup_{0 \leq t \leq T} \norm{\Ec_t}_{5r+10}^{5r+10} \right),
\end{split}
\end{equation*}%
where $C = \max\left\{ \GrowthSum{\mu}{r}{2}^2 , \GrowthSum{\sigma}{r}{2}^4 \right\} \GrowthSum{u}{r}{4}$.
Note that $\int_0^T \int_{\flr{s}}^s d\tau ds = T^2 (2N)^{-1}$.
Thus,
\begin{equation*}
\begin{split}
	\norm{\eta_{1,1}}
	\leq \kappa d C \left( 1 + \sup_{0 \leq t \leq T} \norm{\Ec_t}_{5r+10}^{5r+10} \right) T^2 N^{-1}.
\end{split}
\end{equation*}%
We can also apply It{\^o}'s formula to rewrite $\eta_{1,3}$ as
\begin{equation*}
\begin{split}
	\eta_{1,3}
	&= \int_0^T \Exp{ \GenT{\flr{s}} u(s,\Ec_s,\Ec_{\flr{s}}) - \GenT{\flr{s}} u(s,\Ec_{\flr{s}},\Ec_{\flr{s}}) } ds
	\\
	&= \int_0^T \Exp{ \int_{\flr{s}}^s \GenT{\flr{\tau}} [\GenT{\flr{s}} u(s,\cdot,\Ec_{\flr{s}})](\Ec_{\tau},\Ec_{\flr{\tau}}) d\tau } ds.
\end{split}
\end{equation*}%
By \cref{PGcalc_double_generator}, we obtain the same bound for $\eta_{1,3}$ as for $\eta_{1,1}$ (up to adjusting $\kappa$).
We can estimate $\eta_{1,2}$ by
\begin{equation*}
\begin{split}
	\norm{ \eta_{1,2} } &\leq \Exp{ \int_0^T \norm{  \GenO{\flr{s}} u(s,\Ec_{\flr{s}}) - \GenO{s}u(s,\Ec_{\flr{s}}) } ds }
	\\
	&\leq \Exp{ \int_0^T \norm{ \nabla u(s,\Ec_{\flr{s}}) } \norm{ \mu(\flr{s},\Ec_{\flr{s}}) - \mu(s,\Ec_{\flr{s}}) } ds }
	\\
	&\quad+ \Exp{ \int_0^T \frac{1}{2} \norm{ \sigma(\flr{s},\Ec_{\flr{s}})\sigma(\flr{s},\Ec_{\flr{s}})^T - \sigma(s,\Ec_{\flr{s}})\sigma(s,\Ec_{\flr{s}})^T } \norm{ \nabla^2u(s,\Ec_{\flr{s}}) } ds }.
\end{split}
\end{equation*}%
Furthermore, using the (H{\"o}lder) growth measure and the equality $\int_0^T (s-\flr{s})^{\alpha} ds = (1+\alpha)^{-1} T^{\alpha+1} N^{-\alpha}$,
\begin{equation*}
\begin{split}
	&\Exp{ \int_0^T \norm{ \nabla u(s,\Ec_{\flr{s}}) } \norm{ \mu(\flr{s},\Ec_{\flr{s}}) - \mu(s,\Ec_{\flr{s}}) } ds }
	\\
	&\leq \kappa \GrowthSum{\nabla u}{r+3}{0} \GrowthHoelder{\mu}{r}{\alpha} \left( 1 + \sup_{0 \leq t \leq T} \norm{\Ec_t}_{2r+3}^{2r+3} \right) T^{\alpha+1} N^{-\alpha}.
\end{split}
\end{equation*}%
By \cref{PGcalc_int}, $\GrowthSum{\nabla u}{r+3}{0} \leq \GrowthSum{u}{r}{4}$.
Similarly,
\begin{equation*}
\begin{split}
	&\Exp{ \int_0^T \frac{1}{2} \norm{ \sigma(\flr{s},\Ec_{\flr{s}})\sigma(\flr{s},\Ec_{\flr{s}})^T - \sigma(s,\Ec_{\flr{s}})\sigma(s,\Ec_{\flr{s}})^T } \norm{ \nabla^2u(s,\Ec_{\flr{s}}) } ds }
	\\
	&\leq \frac{1}{2} \Exp{ \int_0^T \big( \norm{ \sigma(\flr{s},\Ec_{\flr{s}}) } + \norm{ \sigma(s,\Ec_{\flr{s}}) } \big) \norm{ \sigma(\flr{s},\Ec_{\flr{s}}) - \sigma(s,\Ec_{\flr{s}}) } \norm{ \nabla^2u(s,\Ec_{\flr{s}}) } ds }
	\\
	&\leq \kappa \GrowthSum{u}{r}{4} \GrowthSum{\sigma}{r}{2} \GrowthHoelder{\sigma}{r}{\alpha} \left( 1 + \sup_{0 \leq t \leq T} \norm{\Ec_t}_{3r+4}^{3r+4} \right) T^{\alpha+1} N^{-\alpha}.
\end{split}
\end{equation*}%
Therefore,
\begin{equation*}
	\norm{ \eta_{1,2} }
	\leq \kappa C' \left( 1 + \sup_{0 \leq t \leq T} \norm{\Ec_t}_{3r+4}^{3r+4} \right) T^{\alpha+1} N^{-\alpha},
\end{equation*}%
where $C' = \max\left\{ \GrowthHoelder{\mu}{r}{\alpha} , \GrowthSum{\sigma}{r}{2} \GrowthHoelder{\sigma}{r}{\alpha} \right\} \GrowthSum{u}{r}{4}$.
Combining the bounds for $\eta_{1,1}$, $\eta_{1,2}$, and $\eta_{1,3}$, we find the following bound for $\eta_1$;
\begin{equation*}
	\norm{\eta_1}
	\leq \kappa d \max\{C,C'\} \left( 1 + \sup_{0 \leq t \leq T} \norm{\Ec_t}_{5r+10}^{5r+10} \right) T^2 N^{-\alpha}.
\end{equation*}%
Next, we bound $\eta_2$, which we can write as $\eta_2 = \sum_{n=0}^{N-1} \int_{t_n}^{t_{n+1}} \Exp{ g(t_n,\Ec_{t_n}) - g(s,\Ec_s) } ds$.
By It{\^o}'s formula,
\begin{equation*}
	\Exp{ g(t_n,\Ec_{t_n}) - g(s,\Ec_s) } = - \Exp{ \int_{t_n}^s \partial_t g(\tau,\Ec_{\tau}) + \GenT{\flr{\tau}} g(\tau,\Ec_{\tau},\Ec_{\flr{\tau}}) d\tau }.
\end{equation*}%
By \cref{PGcalc_int,PGcalc_double_generator},
\begin{equation*}
\begin{split}
	\Exp{ \norm{ \GenT{\flr{\tau}} g(\tau,\Ec_{\tau},\Ec_{\flr{\tau}}) } }
	&\leq \kappa \max\left\{ \GrowthSum{\mu}{r}{2} , \GrowthSum{\sigma}{r}{2}^2 \right\} \GrowthSum{g}{r}{2} \left( 1 + \sup_{0 \leq t \leq T} \norm{\Ec_t}_{3 r + 4}^{3 r + 4} \right).
\end{split}
\end{equation*}%
Thus,
\begin{equation*}
\begin{split}
	\norm{ \eta_2 }
	&\leq \kappa C'' \left( 1 + \sup_{0 \leq t \leq T} \norm{\Ec_t}_{3 r + 4}^{3 r + 4} \right) \sum_{n=0}^{N-1} \int_{t_n}^{t_{n+1}} \int_{t_n}^s d\tau ds
	\leq \kappa C'' \left( 1 + \sup_{0 \leq t \leq T} \norm{\Ec_t}_{3 r + 4}^{3 r + 4} \right) T^2 N^{-1},
\end{split}
\end{equation*}%
where $C'' = \max\left\{ \GrowthSum{\partial_t g}{r}{0} , \max\left\{ \GrowthSum{\mu}{r}{2} , \GrowthSum{\sigma}{r}{2}^2 \right\} \GrowthSum{g}{r}{2} \right\}$.
Combining the bounds for $\eta_1$ and $\eta_2$, we have shown
\begin{equation*}
\begin{split}
	\norm{v(x) - u(0,x)}
	&\leq \norm{\eta_1} + \norm{\eta_2}
	\leq \kappa d C''' \left( 1 + \sup_{0 \leq t \leq T} \norm{\Ec_t}_{5r+10}^{5r+10} \right) T^2 N^{-\alpha},
\end{split}
\end{equation*}%
where $C''' = \max\{C,C',C''\} = \Cc_{r,\alpha}(\mu,\sigma,g,u)$.
By \cref{Euler_scheme_well_defined}, if $N \geq 16 c T$, then
\begin{equation*}
	\sup_{0 \leq t \leq T} \norm{\Ec_t}_{5r+10}^{5r+10}
	\leq e^{\kappa c T} (\kappa c T+\norm{x}^{5r+10}).
\end{equation*}%
This yields the final bound
\begin{equation*}
\begin{split}
	\norm{v(x) - u(0,x)}
	&\leq \kappa d C_{r,\alpha} (1+c) T^3 e^{\kappa c T} \left( 1 + \norm{x}^{5r+10} \right) N^{-\alpha}.
\end{split}
\end{equation*}%
\end{proof}

In \cref{weak_error_rate_thrm}, we see that the error rate is free of the curse of dimensionality as long as we can reasonably control the constants $c$ and $C_{r,\alpha}$.
On the one hand, $C_{r,\alpha}$ appears as a multiplicative factor in the bound and, hence, should grow at most polynomially in the dimension.
On the other hand, $c$ appears as an exponent and, hence, $\GrowthSingle{\mu}$ and $\GrowthSingle{\sigma}$ should grow at most logarithmically in the dimension.
The fact that we have to control these constants is intuitive:
if any of the involved functions itself grows exponentially in the dimension, then we cannot expect the discrete Euler scheme to explore all relevant regions with a time grid whose fineness is polynomial in the dimension.


\subsection{Monte Carlo estimates}

To obtain a numerical scheme from the Feynman--Kac formula, we approximate the expected value with a Monte Carlo method.
The Monte Carlo method is famously known to break the curse of dimensionality.
However, this folklore suppresses the variance of the integrand appearing in the Monte Carlo estimate.
Therein can appear an implicit dependence on the dimension and, hence, may still incur the curse.
We present an estimate of the Monte Carlo method based on the polynomial growth measure to quantify the implicit dependence.

\begin{proposition}
\label{MC_error}
	Let $\Ec_1(x),\dots,\Ec_M(x)$ be $M \in \N$ independent copies of $\Ec^N(x)$.
Consider $v \colon \R^d \rightarrow \R^{\dout}$ and $\phi \colon \Omega \times \R^d \rightarrow \R^{\dout}$ given by
\begin{equation*}
\begin{split}
	v(x) &= \Exp{ f\left(\Ec^N(x)_T\right) + \frac{T}{N} \sum_{n=0}^{N-1} g\left(t_n,\Ec^N(x)_{t_n}\right) },
	\\
	\phi(x) &= \frac{1}{M} \sum_{m=1}^M \left[ f\left(\Ec_m(x)_T\right) + \frac{T}{N} \sum_{n=0}^{N-1} g\left(t_n,\Ec_m(x)_{t_n}\right) \right].
\end{split}
\end{equation*}%
Let $\delta,\epsilon \in (0,1)$, $r \in [1,\infty)$, let $K \subseteq \R^d$ be a compact set, and denote $c = \max\{\GrowthSingle{\mu} , \GrowthSingle{\sigma}^2\}$ and $C_r = \max\{\GrowthSum{f}{r}{0} , \GrowthSum{g}{r}{0}\}$.
If $N \geq 16cT$, then there exists a constant $\kappa_r > 0$, depending only on $r$, such that
\begin{equation*}
\begin{split}
	\Pb\left( \norm{v - \phi}_{L^2(K)} \leq \epsilon \right)
	&\geq 1 - \delta
\end{split}
\end{equation*}%
provided $M$ is so large that
\begin{equation*}
	\epsilon \sqrt{\delta M} \geq \kappa_r C_r \sqrt{(1 + c) T^3} e^{\kappa_r c T} \norm{1+\norm{x}^r}_{L^2(K)}.
\end{equation*}%
\end{proposition}

\begin{proof}
	We use a standard argument for bounding Monte Carlo errors.
First, by Markov's inequality,
\begin{equation*}
\begin{split}
	\Pb\left( \norm{v - \phi}_{L^2(K)} \leq \epsilon \right)
	&\geq 1 - \epsilon^{-2} \int_K \Exp{ \norm{v(x) - \phi(x)}^2 } dx.
\end{split}
\end{equation*}%
Secondly,
\begin{equation*}
\begin{split}
	\Exp{ \norm{v(x) - \phi(x)}^2 }
	&\leq \frac{1}{M} \Exp{ \norm{ f\left(\Ec_1(x)_T\right) + \frac{T}{N} \sum_{n=0}^{N-1} g\left(t_n,\Ec_1(x)_{t_n}\right) }^2 }.
\end{split}
\end{equation*}%
By Jensen's inequality,
\begin{equation*}
	\Exp{ \norm{v(x) - \phi(x)}^2 }
	\leq \frac{2}{MN} \sum_{n=0}^{N-1} \Exp{ \norm{ f\left(\Ec_1(x)_T\right) }^2 + T^2 \norm{ g\left(t_n,\Ec_1(x)_{t_n}\right) }^2 }.
\end{equation*}%
As in the previous proofs, $\kappa \geq 1$ denotes a constant that depends only on $r$ but may change form line to line.
We use the bound $\Exp{ \norm{ g\left(t_n,\Ec_1(x)_{t_n}\right) }^2 } \leq \kappa C_r^2 ( 1 + \norm{\Ec_1(x)_{t_n}}_{2r}^{2r} )$ and similarly for $f(\Ec_1(x)_T)$ to find
\begin{equation*}
	\Exp{ \norm{v(x) - \phi(x)}^2 }
	\leq \frac{\kappa C_r^2 T^2}{M} \left( 1 + \sup_{0 \leq t \leq T} \norm{\Ec_1(x)_t}_{2r}^{2r} \right).
\end{equation*}%
By \cref{Euler_scheme_well_defined}, if $N \geq 16cT$, then
\begin{equation*}
	\sup_{0 \leq t \leq T} \norm{\Ec_1(x)_t}_{2r}^{2r}
	\leq \kappa T e^{\kappa c T} (c+\norm{x}^{2r}).
\end{equation*}%
Thus,
\begin{equation*}
	\int_K \Exp{ \norm{v(x) - \phi(x)}^2 } dx
	\leq \frac{\kappa C_r^2 T^3}{M} (1 + c) e^{\kappa c T} \norm{1+\norm{x}^r}^2_{L^2(K)}.
\end{equation*}%
\end{proof}

By combining the error analysis for the Euler scheme and the Monte Carlo method, we find the following rate of convergence for the Monte Carlo Euler scheme.

\begin{corollary}
\label{MCES_error}
	Let \cref{assumption_weak_error} hold and let $\Ec_1(x),\dots,\Ec_M(x)$ be $M \in \N$ independent copies of $\Ec^N(x)$.
Consider $\phi \colon \Omega \times \R^d \rightarrow \R^{\dout}$ given by
\begin{equation*}
\begin{split}
	\phi(x) &= \frac{1}{M} \sum_{m=1}^M \left[ f\left(\Ec_m(x)_T\right) + \frac{T}{N} \sum_{n=0}^{N-1} g\left(t_n,\Ec_m(x)_{t_n}\right) \right].
\end{split}
\end{equation*}%
Let $\alpha \in (0,1]$, $\delta,\epsilon \in (0,1)$, $r \in [0,\infty)$, let $K \subseteq \R^d$ be compact, and denote $c = \max\{\GrowthSingle{\mu},\GrowthSingle{\sigma}^2\}$ and $c_{r,\alpha} = \max\{ \GrowthSum{f}{r}{0} , \GrowthSum{g}{r}{0} , \Cc_{r,\alpha}(\mu,\sigma,g,u) \}$; see \eqref{cumbersome_constant}.
If $N \geq 16cT$, then there exists a constant $\kappa_r > 0$, depending only on $r$, such that
\begin{equation*}
\begin{split}
	\Pb\left( \norm{u(0,\cdot) - \phi}_{L^2(K)} \leq \epsilon \right)
	&\geq 1 - \delta
\end{split}
\end{equation*}%
provided $M$ and $N$ are so large that
\begin{equation*}
	\frac{\epsilon \sqrt{\delta}}{M^{-1/2} + \sqrt{\delta} N^{-\alpha}}
	\geq d c_{r,\alpha} \kappa_r (1+c) T^3 e^{ \kappa_r c T} \norm{1 + \norm{x}^{5r+10}}_{L^2(K)}.
\end{equation*}%
\end{corollary}

\begin{proof}
	This is immediate from \cref{weak_error_rate_thrm,MC_error}.
\end{proof}

As before, we see that the error rate for the Monte Carlo Euler scheme is free of the curse of dimensionality as long as we can control $c$ and $c_{r,\alpha}$.
The same comment we made after the proof of \cref{weak_error_rate_thrm} applies here.
An implementation of this scheme requires $M N$ independent simulations of a standard $d$-dimensional normal random vector.

\begin{remark}
\label{MCES_time_dependent}
	If we consider the Euler scheme starting at time $t \in [0,T]$ instead of at time 0, that is setting $t_n = t + n(T-t)/N$ and defining
\begin{equation*}
\begin{split}
	\Ec^N(x)_s &= \Ec^N(x)_{t_n} + \mu(t_n,\Ec^N(x)_{t_n}) (s-t_n) + \sigma(t_n,\Ec^N(x)_{t_n}) (W_s-W_{t_n}), \qquad t_n \leq s \leq t_{n+1},
	\\
	\Ec^N(x)_{t} &= x,
\end{split}
\end{equation*}%
then the Monte Carlo Euler scheme in \cref{MCES_error} becomes a time-dependent approximation $\phi(t,x)$ with
\begin{equation*}
\begin{split}
	\Pb\left( \norm{u - \phi}_{L^2([0,T] \times K)} \leq \epsilon \right)
	&\geq 1 - \delta
\end{split}
\end{equation*}%
provided $M$ and $N$ are so large that
\begin{equation*}
	\frac{\epsilon \sqrt{\delta}}{M^{-1/2} + \sqrt{\delta} N^{-\alpha}}
	\geq d c_{r,\alpha} \kappa_r (1+c) T^4 e^{ \kappa_r c T} \norm{1 + \norm{x}^{5r+10}}_{L^2(K)}.
\end{equation*}%
To implement the time-dependent scheme, $M N$ independent simulations of a $d$-dimensional normal random vector are still sufficient because we do not need the simulations of $W_{t_{n+1}}-W_{t_n} \sim \Nc(0,(T-t)/N \cdot \Ind_{\dWiener})$ to be independent for different $t$, so we can use simulations of $\sqrt{(T-t)/N} Z$ with the same $Z \sim \Nc(0,\Ind_{\dWiener})$ for each $t$.
\end{remark}

This concludes the discussion about the numerical scheme for approximating PDE \eqref{Intro_PDE} in $L^2$.
In the next section, we will set the stage to obtain Sobolev approximations of the PDE.


\section{Numerical Sobolev approximation of PDE solutions}
\label{section_comDer}

\subsection{Augmented derivatives of PDEs}

In the introduction, we motivated considering the combined process $(X,DX)$, where $X$ solves SDE \eqref{Intro_SDE}.
The combined process solves another SDE, namely $dZ_t = \ED{\mu}(t,Z_t) dt + \ED{\sigma}(t,Z_t) dW_t$ with the drift coefficient $\ED{\mu}(t,x,y) = (\mu(t,x),\nabla\mu(t,x)^Ty)$ and diffusion coefficient $\ED{\sigma}(t,x,y) = (\sigma(t,x),\nabla\sigma(t,x)^Ty)$.
In light of this, it becomes natural to introduce the following terminology of an augmented derivative.
This terminology is chiefly a matter of notation, which will enable us to rewrite the problem of obtaining Sobolev approximations of PDE \eqref{Intro_PDE} so as to fit the setting of the previous section.

\begin{definition}[Augmented derivative]
	Given a Fr{\'e}chet differentiable function $\fct \colon \Bc_1 \rightarrow \Bc_2$ between two Banach spaces, whose Fr{\'e}chet derivative at a point $x_1 \in \Bc_1$ is denoted $\Bc_1 \rightarrow \Bc_2$, $x_2 \mapsto D\fct_{x_1}(x_2)$, we let $\ED{\fct} \colon \Bc_1 \times \Bc_1 \rightarrow \Bc_2 \times \Bc_2$ be the function
\begin{equation*}
	\ED{\fct}(x_1,x_2) =
	\begin{pmatrix}
	\fct(x_1) \\
	D\fct_{x_1}(x_2)
	\end{pmatrix}.
\end{equation*}%
If the function $\fct \colon [0,T] \times \Bc_1 \rightarrow \Bc_2$ is time-dependent, then we denote $\ED{\fct}(t,x) = \WED{\fct(t,\cdot)}(x)$.
\end{definition}

The augmented derivative augments the derivative $D\fct$ with the original function.
Note that it is linear in $h$ and satisfies the chain rule $\WED{\fct_1 \circ \fct_2} = \ED{\fct}_1 \circ \ED{\fct}_2$.
The convenience of the augmented derivative in terms of notation is demonstrated in the next result.
In general, the standard derivative does {\it not} satisfy $\nabla L_{\mu,\sigma}\fct = L_{\nabla\mu,\nabla\sigma}\nabla \fct$.
This equality becomes true when we replace the standard derivative by the augmented derivative.
We briefly omit the time-dependence of $\mu$ and $\sigma$ for increased readability.

\begin{proposition}[Commutativity of augmented derivatives and generators]
\label{comDer_generators_commute}
	Suppose $\mu$ and $\sigma$ are differentiable in the space variable and $\fct \colon \R^d \rightarrow \R^{\dout}$ is three times differentiable.
Let $\tau \colon \R^{4d} \rightarrow \R^{4d}$ be given by $\tau(x_1,x_2,x_3,x_4) = (x_1,x_3,x_2,x_4)$.
Then $\WED{\GenT{} \fct} = \GenDiffT{} \ED{\fct} \circ \tau$ and $\WED{\GenO{} \fct} = \GenDiffO{} \ED{\fct}$.
\end{proposition}

\begin{proof}
	On the one hand, we calculate
\begin{equation*}
	\nabla\ED{\fct}(x)^T\ED{\mu}(y) =
	\begin{pmatrix}
	\nabla \fct(x_1)^T \mu(y_1) \\
	[\nabla^2\fct(x_1)^Tx_2]\mu(y_1) + \nabla \fct(x_1)^T \nabla \mu(y_1)^T y_2
	\end{pmatrix}
\end{equation*}%
and
\begin{equation*}
	\tr\big( \ED{\sigma}(y) \ED{\sigma}(y)^T \nabla^2 \ED{\fct}(x) \big) =
	\begin{pmatrix}
	\tr\big( \sigma(y_1) \sigma(y_1)^T \nabla^2\fct(x_1) \big)
	\\
	\tr\big( \sigma(y_1) \sigma(y_1)^T \nabla^3 \fct(x_1)^T x_2 \big) + 2 \tr\big( [\nabla\sigma(y_1)^Ty_2] \sigma(y_1)^T \nabla^2 \fct(x_1) \big)
	\end{pmatrix}.
\end{equation*}%
On the other hand,
\begin{equation*}
\begin{split}
	\nabla \GenT{}\fct(x_1,y_1)^T \begin{pmatrix} x_2 \\ y_2 \end{pmatrix}
	&=
	[\nabla^2\fct(x_1)^Tx_2]\mu(y_1) + \frac{1}{2} \tr\big( \sigma(y_1) \sigma(y_1)^T \nabla^3 \fct(x_1)^T x_2 \big)
	\\
	&\quad+ \nabla \fct(x_1)^T \nabla \mu(y_1)^T y_2 + \tr\big( [\nabla\sigma(y_1)^Ty_2] \sigma(y_1)^T \nabla^2 \fct(x_1) \big).
\end{split}
\end{equation*}%
Thus, $\WED{\GenT{} \fct} = \GenDiffT{} \ED{\fct} \circ \tau$.
Lastly, let $\iota_d \colon \R^d \rightarrow \R^{2d}$ and $\iota_{2d} \colon \R^{2d} \rightarrow \R^{4d}$ be given by $\iota_d(x) = (x,x)$ and $\iota_{2d}(x,y) = (x,y,x,y)$, respectively.
Then $\GenO{}\fct = \GenT{}\fct \circ \iota_d$ and $\tau \circ \ED{\iota}_d = \iota_{2d}$.
In particular, by the chain rule for the augmented derivative,
\begin{equation*}
	\WED{\GenO{} \fct} = \WED{\GenT{} \fct} \circ \ED{\iota}_d = \GenDiffT{} \ED{\fct} \circ \tau \circ \ED{\iota}_d = \GenDiffT{} \ED{\fct} \circ \iota_{2d} = \GenDiffO{} \ED{\fct}.
\end{equation*}%
\end{proof}

As an immediate consequence of the previous result, we find that the augmented derivative of a solution to PDE \eqref{Intro_PDE} solves the augmented derivative of the PDE, meaning the PDE whose coefficients are the augmented derivatives of the original coefficients.

\begin{corollary}[Commutativity of augmented derivatives and PDEs]
\label{comDer_solves_comPDE}
	Suppose $\mu$, $\sigma$, and $g$ are differentiable in the space variable, and $u \in C^{1,3}([0,T] \times \R^d,\R^{\dout})$ solves the PDE $\partial_t u + \GenO{}u + g = 0$.
Then, $\ED{u}$ solves the PDE $\partial_t \ED{u} + \GenDiffO{}\ED{u} + \ED{g} = 0$.
\end{corollary}

\begin{proof}
	This is an immediate consequence of \cref{comDer_generators_commute}.
\end{proof}

In the next section, we will see that the augmented derivative interacts as nicely with SDEs as with PDEs.


\subsection{Augmented derivatives of Euler schemes}

Consider the SDE
\begin{equation}
\label{SDE}
\begin{split}
	X_t &= x + \int_0^t \mu(s,X_s) ds + \int_0^t \sigma(s,X_s) dW_s, \qquad 0 \leq t \leq T,
	\\
	X_0 &= x,
\end{split}
\end{equation}%
with initial condition $x \in \R^d$.
It is a classical result that this SDE admits a unique solution $\SDE(x) \in \Einf{}$ under certain regularity assumptions on $\mu$ and $\sigma$.
Under stronger regularity assumptions on the coefficients, the solution map $\R^d \rightarrow \CAd{d}{p}$, $x \mapsto \SDE(x)$ becomes differentiable;
\cite{ImkeDosReisSalk2019}.
In particular, we can take its augmented derivative.
The augmented derivative of a Fr{\'e}chet differentiable map $\R^d \rightarrow \CAd{d}{p}$ is a map $\R^{2d} \rightarrow \CAd{d}{p} \times \CAd{d}{p}$.
We can identify the product $\CAd{d}{p} \times \CAd{d}{p}$ with the space $\CAd{2d}{p}$.
This identification works on the the level of Banach spaces because the norm $\normCAd{X}{p} + \normCAd{Y}{p}$ on $(X,Y) \in \CAd{d}{p} \times \CAd{d}{p}$ is equivalent to the norm $\normCAd{(X,Y)}{p}$ on $(X,Y) \in \CAd{2d}{p}$.
Henceforth, we regard the augmented derivative of a map $\R^d \rightarrow \CAd{d}{p}$ as taking values in $\CAd{2d}{p}$.

\begin{assumption}
\label{assumption_SDE_diff}
	Suppose $\mu$ and $\sigma$ are continuous in time and continuously differentiable in space such that their spatial derivatives are globally bounded uniformly in time.
\end{assumption}

\begin{proposition}[Commutativity of augmented derivatives and SDEs]
\label{comDer_SDE_solutions_commute}
	Let \cref{assumption_SDE_diff} hold.
Then, for any $p \in [2,\infty)$, the function $\SDE \colon \R^d \rightarrow \CAd{d}{p}$ is Fr{\'e}chet differentiable and $\WED{\SDE} = \SDEED$.
In particular, $\WED{\SDE}$ maps into $\Einf{2}$ and is independent of $p$.
\end{proposition}

\begin{proof}
	This is proved in \cite{ImkeDosReisSalk2019}, only rephrased in the language of augmented derivatives.
\end{proof}

Naturally, the same assumptions on the coefficients that guaranteed differentiability of $x \mapsto \SDE(x)$ provide the same regularity for the Euler scheme.
Again, this is conveniently formulated with the help of augmented derivatives.
For a moment, we make the dependence of the Euler scheme on the coefficient functions explicit and write $\Euler = \Ec^N$.

\begin{proposition}[Commutativity of augmented derivatives and Euler schemes]
\label{comDer_Euler_commute}
	Let \cref{assumption_SDE_diff} hold.
Then, for any $p \in [2,\infty)$, the function $\Euler \colon \R^d \rightarrow \CAd{d}{p}$ is Fr{\'e}chet differentiable and $\WED{\Euler} = \EulED$.
In particular, $\WED{\Euler}$ maps into $\Einf{2}$ and is independent of $p$.
\end{proposition}

\begin{proof}
	By \cref{assumption_SDE_diff}, the functions $\mu$ and $\sigma$ are Lipschitz continuous with a Lipschitz constant independent of time.
This implies that $\ED{\mu}$ and $\ED{\sigma}$ grow at most linearly in the space variable.
In particular, \cref{Euler_scheme_well_defined} and its proof assert that the Euler scheme $\EulED$ for the functions $\ED{\mu}$ and $\ED{\sigma}$ is indeed an element of $\Einf{2}$.
For any $x,y \in \R^d$, let $D\Ec(x,y)$ be the second component of $\EulED$, that is
\begin{equation*}
	D\Ec(x,y)_t
	= y + \int_0^t \nabla \mu(\flr{s},\Euler(x)_{\flr{s}})^T D\Ec(x,y)_{\flr{s}} ds + \int_0^t \nabla \sigma(\flr{s},\Euler(x)_{\flr{s}})^T D\Ec(x,y)_{\flr{s}} dW_s.
\end{equation*}%
Fix $p \in [2,\infty)$.
We show that $D\Ec(x,y)$ is the Fr{\'e}chet derivative of $\Ec = \Euler \colon \R^d \rightarrow \CAd{d}{p}$ at $x$ applied to $y$.
In this proof, there is no need to track constants explicitly.
Henceforth, $c$ denotes a constant that may change from line to line but always stays independent of $x$ and $y$.
Abbreviate
\begin{equation*}
	R_n = \Exp{ \sup_{0 \leq t \leq t_n} \norm{\Ec(x+y)_t-\Ec(x)_t-D\Ec(x,y)_t}^p }.
\end{equation*}%
Our goal is to show that $\lim_{y \rightarrow 0} R_N^{1/p}\norm{y}^{-1} = 0$.
Using the recursive structure of the Euler scheme, we see that
\begin{equation*}
\begin{split}
	R_{n+1}
	&\leq c R_n
	+ c \Exp{ \norm{ \mu(t_n,\Ec(x+y)_{t_n}) - \mu(t_n,\Ec(x)_{t_n}) - \nabla \mu(t_n,\Ec(x)_{t_n})^T D\Ec(x,y)_{t_n} }^p }
	\\
	&\quad+ c \Exp{ \sup_{t_n \leq t \leq t_{n+1}} \norm{ \left( \sigma(t_n,\Ec(x+y)_{t_n}) - \sigma(t_n,\Ec(x)_{t_n}) - \nabla \sigma(t_n,\Ec(x)_{t_n})^T D\Ec(x,y)_{t_n}\right) (W_t-W_{t_n}) }^p }.
\end{split}
\end{equation*}%
To the last term, we can apply the Burkholder--Davis--Gundy inequality to find
\begin{equation*}
\begin{split}
	&\Exp{ \sup_{t_n \leq t \leq t_{n+1}} \norm{ \left( \sigma(t_n,\Ec(x+y)_{t_n}) - \sigma(t_n,\Ec(x)_{t_n}) - \nabla \sigma(t_n,\Ec(x)_{t_n})^T D\Ec(x,y)_{t_n} \right) (W_t-W_{t_n}) }^p }
	\\
	&\leq c \Exp{ \sup_{t_n \leq t \leq t_{n+1}} \norm{ \sigma(t_n,\Ec(x+y)_{t_n}) - \sigma(t_n,\Ec(x)_{t_n}) - \nabla \sigma(t_n,\Ec(x)_{t_n})^T D\Ec(x,y)_{t_n} }^p }.
\end{split}
\end{equation*}%
Let us write
\begin{equation*}
	\eta_n(v,w) = \sup_{r \in [0,1]} \norm{\nabla\mu(t_n,(1-r)v+rw)-\nabla\mu(t_n,v)} + \sup_{r \in [0,1]} \norm{\nabla\sigma(t_n,(1-r)v+rw)-\nabla\sigma(t_n,v)}
\end{equation*}%
so that uniformly
\begin{equation*}
\begin{split}
	\norm{\mu(t_n,w)-\mu(t_n,v)-\nabla\mu(t_n,v)^T(w-v)} + \norm{\sigma(t_n,w)-\sigma(t_n,v)-\nabla\sigma(t_n,v)^T(w-v)}
	&\leq \eta_n(v,w) \norm{v-w}.
\end{split}
\end{equation*}%
Using that $\mu$ is Lipschitz continuous, we estimate
\begin{equation*}
\begin{split}
	&\norm{ \mu(t_n,\Ec(x+y)_{t_n}) - \mu(t_n,\Ec(x)_{t_n}) - \nabla \mu(t_n,\Ec(x)_{t_n})^T D\Ec(x,y)_{t_n} }^p \\
	&\leq c \norm{ \nabla \mu(t_n,\Ec(x)_{t_n})^T (\Ec(x+y)_{t_n}-\Ec(x)_{t_n} - D\Ec(x,y)_{t_n}) }^p
	\\
	&\quad + c \norm{ \mu(t_n,\Ec(x+y)_{t_n}) - \mu(t_n,\Ec(x)_{t_n}) - \nabla \mu(t_n,\Ec(x)_{t_n})^T (\Ec(x+y)_{t_n}-\Ec(x)_{t_n}) }^p
	\\
	&\leq c \norm{\Ec(x+y)_{t_n}-\Ec(x)_{t_n} - D\Ec(x,y)_{t_n} }^p + c \eta_n(\Ec(x)_{t_n},\Ec(x+y)_{t_n})^p \norm{ \Ec(x+y)_{t_n}-\Ec(x)_{t_n} }^p.
\end{split}
\end{equation*}%
We proceed analogously for the term involving $\sigma$.
Thus,
\begin{equation}
\label{Pf_comDer_Euler_commute_recursion}
\begin{split}
	R_{n+1}
	&\leq c R_n + c \Exp{ \eta_n(\Ec(x)_{t_n},\Ec(x+y)_{t_n})^p \norm{ \Ec(x+y)_{t_n}-\Ec(x)_{t_n} }^p }
	\\
	&\leq c R_n + c \Exp{ \eta_n(\Ec(x)_{t_n},\Ec(x+y)_{t_n})^{2p} }^{1/2} \Exp{ \norm{ \Ec(x+y)_{t_n}-\Ec(x)_{t_n} }^{2p} }^{1/2}.
\end{split}
\end{equation}%
where we used H{\"o}lder's inequality in the second line.
Now, abbreviate
\begin{equation*}
	S_n = \Exp{ \sup_{0 \leq t \leq t_n} \norm{ \Ec(x+y)_t-\Ec(x)_t }^{2p} }.
\end{equation*}%
Similarly as for $R_n$, by the Burkholder-Davis-Gundy inequality and Lipschitz continuity, $S_{n+1} \leq c S_n$.
In particular, $S_n \leq c S_0 = c \norm{y}^{2p}$.
Iterating the recursion in \eqref{Pf_comDer_Euler_commute_recursion}, we find
\begin{equation*}
\begin{split}
	R_N
	&\leq c \norm{y}^p \sum_{n=0}^{N-1} \Exp{ \eta_n(\Ec(x)_{t_n},\Ec(x+y)_{t_n})^{2p} }^{1/2}.
\end{split}
\end{equation*}%
The limit of $\eta_n(\Ec(x)_{t_n},\Ec(x+y)_{t_n})$ as $y \rightarrow 0$ exists pointwise and equals 0.
Furthermore, the Lipschitz continuity of $\mu$ and $\sigma$ implies that $\eta_n$ is uniformly bounded.
Hence, we can apply dominated convergence to find
\begin{equation*}
\begin{split}
	\lim_{y \rightarrow 0} \frac{R_N^{1/p}}{\norm{y}}
	&\leq \lim_{y \rightarrow 0} c \left( \sum_{n=0}^{N-1} \Exp{ \eta_n(\Ec(x)_{t_n},\Ec(x+y)_{t_n})^{2p} }^{1/2} \right)^{1/p}
	= 0.
\end{split}
\end{equation*}%
This concludes Fr{\'e}chet differentiability of $\Ec$ with the correct derivative if we show that $D\Ec(x,y)$ is a bounded linear operator in $y$.
Boundedness was shown implicitly in the previous calculations.
Indeed, by Jensen's inequality,
\begin{equation*}
\begin{split}
	\frac{\normCAd{ D\Ec(x,y) }{p}^p}{\norm{y}^p}
	&\leq c \frac{\sqrt{S_N}+R_N}{\norm{y}^p}
	\leq c + c \sum_{n=0}^{N-1} \Exp{ \eta_n(\Ec(x)_{t_n},\Ec(x+y)_{t_n})^{2p} }^{1/2},
\end{split}
\end{equation*}%
and this is bounded uniformly in $y$ because $\eta_n$ is.
For linearity, abbreviate
\begin{equation*}
	L_n = \Exp{ \sup_{0 \leq t \leq t_n} \norm{D\Ec(x,y+\lambda z)_t-D\Ec(x,y)_t-\lambda D\Ec(x,z)_t}^p }.
\end{equation*}%
As for $S_n$, by the Burkholder--Davis--Gundy inequality and the Lipschitz continuity of $\mu$ and $\sigma$, we obtain $L_N \leq c L_0 = 0$.
\end{proof}

The final preparation for the Sobolev approximation is a standard argument for interchanging expectation and differentiation.

\begin{lemma}
\label{existence_of_compound_derivative}
	Let \cref{assumption_SDE_diff} hold and suppose $f$ and $g$ are continuously differentiable in space with $\GrowthSum{f}{r}{1},\GrowthSum{g}{r}{1} < \infty$ for some $r \geq 0$.
Then, the function $v \colon \R^d \rightarrow \R^{\dout}$ given by
\begin{equation*}
	v(x) = \Exp{ f\left(\Euler(x)_T\right) + \frac{T}{N} \sum_{n=0}^{N-1} g\left(t_n,\Euler(x)_{t_n}\right) }
\end{equation*}%
is differentiable and its augmented derivative is
\begin{equation*}
	\ED{v}(x) = \Exp{\ED{f}\left(\EulED(x)_T\right)  + \frac{T}{N} \sum_{n=0}^{N-1} \ED{g}\left(t_n,\EulED(x)_{t_n}\right) }.
\end{equation*}%
\end{lemma}

\begin{proof}
	Let us write $D\Ec$ for the Fr{\'e}chet derivative of $\Ec = \Euler$ and
\begin{equation*}
	\eta(v,w) = \sup_{r \in [0,1]} \norm{\nabla f((1-r)v+rw)-\nabla f(v)}.
\end{equation*}%
We will show differentiability of the map $x \mapsto F(x) = \Exp{f(\Ec(x)_T)}$.
The argument for $x \mapsto \Exp{g(t_n,\Ec(x)_{t_n})}$ is completely analogous.
We begin with the estimate
\begin{equation*}
\begin{split}
	\epsilon(x,y)
	&= \norm{ F(x+y) - F(x) - \Exp{ \nabla f(\Ec(x)_T)^T D\Ec(x,y)_T }}
	\\
	&\leq \Exp{\norm{ \nabla f(\Ec(x)_T)^T(\Ec(x+y)_T-\Ec(x)_T-D\Ec(x,y)_T) }}
	\\
	&\quad + \Exp{\norm{ f(\Ec(x+y)_T) - f(\Ec(x)_T) - \nabla f(\Ec(x)_T)^T(\Ec(x+y)_T-\Ec(x)_T) }}.
\end{split}
\end{equation*}%
By H{\"o}lder's inequality and the definition of $\eta$,
\begin{equation*}
\begin{split}
	\epsilon(x,y)
	&\leq \Exp{\norm{ \nabla f(\Ec(x)_T) }^2}^{1/2} \Exp{\norm{ \Ec(x+y)_T-\Ec(x)_T-D\Ec(x,y)_T }^2}^{1/2}
	\\
	&\quad + \Exp{ \eta(\Ec(x)_T,\Ec(x+y)_T)^2 }^{1/2} \Exp{ \norm{ \Ec(x+y)_T-\Ec(x)_T }^2}^{1/2}.
\end{split}
\end{equation*}%
The expectation of $\norm{ \nabla f(\Ec(x)_T) }^2$ is finite since $\Ec(x) \in \Einf{}$ and $\GrowthSum{f}{r}{1} < \infty$ for some $r \geq 0$.
As in the proof of \cref{comDer_Euler_commute}, we find $\Exp{ \norm{ \Ec(x+y)_T-\Ec(x)_T }^2} \leq c \norm{y}^2$ for some constant $c$ independent of $x$ and $y$.
Thus, differentiability of $\Ec$ implies
\begin{equation*}
\begin{split}
	\lim_{y \rightarrow 0} \frac{ \epsilon(x,y) }{\norm{y}}
	\leq \lim_{y \rightarrow 0} \sqrt{c}\ \Exp{ \eta(\Ec(x)_T,\Ec(x+y)_T)^2 }^{1/2}.
\end{split}
\end{equation*}%
To this term, we can apply dominated convergence because $\GrowthSum{f}{r}{1} < \infty$ and
\begin{equation*}
\begin{split}
	\Exp{ \eta(\Ec(x)_T,\Ec(x+y)_T)^2 }
	&\leq 4 \GrowthSum{f}{r}{1}^2 \Exp{ ( 1 + \norm{\Ec(x)_T} + \norm{\Ec(x+y)_T-\Ec(x)_T} )^{2r} }
	\\
	&\leq c \GrowthSum{f}{r}{1}^2 \left( 1 + \normCAd{\Ec(x)}{2r}^{2r} + \norm{y}^{2r} \right)
\end{split}
\end{equation*}%
for some (other) constant $c$ independent of $x$ and $y$.
We have seen in the proof of \cref{comDer_Euler_commute} that $D\Ec(x,y)$ and, hence, $\Exp{ \nabla f(\Ec(x)_T)^T D\Ec(x,y)_T }$ is a bounded linear map in $y$.
\end{proof}


\subsection{Sobolev error rate for the Monte Carlo Euler scheme}

In this section, we prove the error rate for the Sobolev approximation of PDE \eqref{Intro_PDE}.
The results of the previous section enable us to directly apply the weak error rate of the Euler scheme and the Monte Carlo estimate to $\ED{u}$ in place of $u$.
Since the augmented derivatives of all functions involved have to satisfy the assumptions of \cref{weak_error_rate_thrm}, we need to upgrade \cref{assumption_weak_error} to one more order of differentiability.

\begin{assumption}
\label{Sobolev_assumption_weak_error}
	Let \cref{assumption_SDE_diff} hold and suppose $\mu$, $\sigma$, and $g$ are three times continuously differentiable in the space variable, $g$ is also continuously differentiable in the time variable, and $f$ is five times continuously differentiable.
Suppose $u \in C^{1,5}([0,T] \times \R^d,\R^{\dout})$ solves the PDE
\begin{equation*}
\begin{split}
	\partial_t u + \GenO{}u + g &= 0 \quad \text{on } [0,T] \times \R^d, \\
	u(T,\cdot) &= f \quad \text{on } \R^d.
\end{split}
\end{equation*}%
\end{assumption}

With these assumptions in place, \cref{weak_error_rate_thrm} can be upgraded immediately to yield an analogous error bound for $\norm{ \ED{u}(0,x) - \ED{v}(x) }$.
Likewise, by applying \cref{MCES_error} to the augmented derivatives of the PDE solution $u$ and its Monte Carlo Euler approximation $\phi$, we find a result of the form
\begin{equation*}
\begin{split}
	\Pb\left( \norm{\ED{u}(0,\cdot) - \ED{\phi}}_{L^2(K)} \leq \epsilon \right)
	&\geq 1 - \delta
\end{split}
\end{equation*}%
for any compact $K \subseteq \R^{2d}$ and sufficiently large $M$ and $N$.
If $K = K_1 \times K_2 \subseteq \R^d \times \R^d$, then the error is given by
\begin{equation*}
\begin{split}
	\norm{\ED{u}(0,\cdot) - \ED{\phi}}_{L^2(K)}^2
	&= \int_{K_1} \left( \norm{u(0,x_1) - \phi(x_1)}^2 + \int_{K_2} \norm{\nabla u(0,x_1)^Tx_2 - \nabla \phi(x_1)^Tx_2}^2 dx_2 \right) dx_1.
\end{split}
\end{equation*}%
In particular, it involves the integral $\int_{K_2} \norm{\nabla u(0,x_1)^Tx_2 - \nabla \phi(x_1)^Tx_2}^2 dx_2$.
This is not how the error of the derivative is usually measured.
Instead, we would like to measure $\norm{\nabla u(0,x_1) - \nabla \phi(x_1)}^2$ in the Frobenius norm.
Converting the former measure of error to the latter can be done as follows.

\begin{lemma}
\label{covert_derivative_error}
	Let $\fct \colon \R^d \rightarrow \R^{\dout}$ be differentiable.
Then, for all $x \in \R^d$,
\begin{equation*}
\begin{split}
	\norm{\nabla \fct(x)}^2
	&\leq \frac{2 d}{\pi^{d/2}} \Gamma\left( \frac{d+4}{2} \right) \int_{\B^d} \norm{\nabla \fct(x)^Ty}^2 dy,
\end{split}
\end{equation*}%
where $\B^d \subseteq \R^d$ is the Euclidean $d$-dimensional ball of radius 1 centered at the origin and $\Gamma$ is the Gamma function $\Gamma(z) = \int_0^{\infty} t^{z-1} e^{-t} dt$.
\end{lemma}

\begin{proof}
	Take $w \in \R^d$ with unit norm such that $\norm{\nabla \fct(x)^Tw}$ equals the operator norm of $\nabla \fct(x)$.
Let $w^{\perp}$ be the orthogonal complement of $w$ in $\R^d$ with respect to the standard inner product.
Then, for any $r \in \R$ and $z \in w^{\perp}$,
\begin{equation*}
\begin{split}
	\norm{\nabla \fct(x)^T(rw+z)}^2 + \norm{\nabla \fct(x)^T(rw-z)}^2
	&\geq 2 r^2 \norm{\nabla \fct(x)^Tw}^2.
\end{split}
\end{equation*}%
Let $w^{\perp}_+ = \{z \in w^{\perp} \colon z^Tw^{\perp}_0 \geq 0\}$ for some fixed $w^{\perp}_0 \in w^{\perp}$.
Then,
\begin{equation*}
\begin{split}
	\int_{\B^d} \norm{\nabla \fct(x)^Ty}^2 dy
	&= \int_{\R} \int_{w^{\perp}} \Ind_{\B^d}(rw+z) \norm{\nabla \fct(x)^T(rw+z)}^2 dz dr
	\\
	&= \int_{\R} \int_{w^{\perp}_+} \Ind_{\B^d}(rw+z) \left( \norm{\nabla \fct(x)^T(rw+z)}^2 + \norm{\nabla \fct(x)^T(rw-z)}^2 \right) dz dr,
\end{split}
\end{equation*}%
where we used that $\Ind_{\B^d}(rw-z) = \Ind_{\B^d}(rw+z)$ by symmetry of the ball $\B^d$.
Thus,
\begin{equation*}
\begin{split}
	\int_{\B^d} \norm{\nabla \fct(x)^Ty}^2 dy
	&\geq 2 \norm{\nabla \fct(x)^Tw}^2 \int_{-1}^1 r^2 \int_{w^{\perp}_+} \Ind_{\B^d}(rw+z) dz dr.
\end{split}
\end{equation*}%
The integral $\int_{w^{\perp}_+} \Ind_{\B^d}(rw+z) dz$ is the volume of half of a  $(d-1)$-dimensional ball of radius $\sqrt{1-r^2}$, so
\begin{equation*}
\begin{split}
	\int_{-1}^1 r^2 \int_{w^{\perp}_+} \Ind_{\B^d}(rw+z) dz dr
	&= \frac{1}{2} \mathrm{Vol}{(\B^{d-1})} \int_{-1}^1 r^2 (1-r^2)^{(d-1)/2} dr
	\\
	&= \frac{1}{2} \mathrm{Vol}{(\B^{d-1})} \frac{\sqrt{\pi}}{2} \frac{\Gamma\left( \frac{d+1}{2} \right)}{\Gamma\left( \frac{d+4}{2} \right)}
	= \frac{\pi^{d/2}}{4} \frac{1}{\Gamma\left( \frac{d+4}{2} \right)}.
\end{split}
\end{equation*}%
Finally, $\norm{\nabla \fct(x)}^2 \leq d \norm{\nabla \fct(x)^Tw}^2$ because the rank of $\nabla \fct(x)$ can be at most $d$.
\end{proof}

Now, we can deduce the main result of this article about the Sobolev approximation of PDE \eqref{Intro_PDE} with a Monte Carlo Euler scheme.

\begin{theorem}
\label{Sobolev_MCES_error}
	Let \cref{Sobolev_assumption_weak_error} hold and let $\Ec_1(x),\dots,\Ec_M(x)$ be $M \in \N$ independent copies of $\Ec^N(x)$.
Consider $\phi \colon \Omega \times \R^d \rightarrow \R^{\dout}$ given by
\begin{equation*}
\begin{split}
	\phi(x) &= \frac{1}{M} \sum_{m=1}^M \left[ f\left(\Ec_m(x)_T\right) + \frac{T}{N} \sum_{n=0}^{N-1} g\left(t_n,\Ec_m(x)_{t_n}\right) \right].
\end{split}
\end{equation*}%
Let $\alpha \in (0,1]$, $\delta,\epsilon \in (0,1)$, $r \in [0,\infty)$, let $K \subseteq \R^d$ be compact, and denote $c = \max\{\GrowthSingle{\ED{\mu}},\GrowthSingle{\ED{\sigma}}^2\}$ and $c_{r,\alpha} = \max\{ \GrowthSum{\ED{f}}{r}{0} , \GrowthSum{\ED{g}}{r}{0} , \Cc_{r,\alpha}(\ED{\mu},\ED{\sigma},\ED{g},\ED{u}) \}$; see \eqref{cumbersome_constant}.
If $N \geq 16cT$, then there exists a constant $\kappa_r > 0$, depending only on $r$, such that
\begin{equation*}
\begin{split}
	\Pb\left( \norm{u(0,\cdot) - \phi}_{L^2(K)} + \norm{\nabla u(0,\cdot) - \nabla\phi}_{L^2(K)} \leq \epsilon \right)
	&\geq 1 - \delta
\end{split}
\end{equation*}%
provided $M$ and $N$ are so large that
\begin{equation*}
	\frac{\epsilon \sqrt{\delta}}{M^{-1/2} + \sqrt{\delta} N^{-\alpha}}
	\geq d^2 c_{r,\alpha} \kappa_r (1+c) T^3 e^{ \kappa_r c T} \norm{1 + \norm{x}^{5r+10}}_{L^2(K)}.
\end{equation*}%
\end{theorem}

\begin{proof}
	Let $\B^d \subseteq \R^d$ be the $d$-dimensional Euclidean ball of radius 1 centered at the origin.
As usual, $\kappa \geq 1$ denotes a constant that depends only on $r$ but may change from line to line.
Observe that $\ED{\Ec}_1(x),\dots,\ED{\Ec}_M(x)$ are independent copies of $\WED{\Ec^N}(x)$.
We apply \cref{MCES_error}, which is possible by
\cref{comDer_solves_comPDE,comDer_Euler_commute,existence_of_compound_derivative},
to find
\begin{equation*}
\begin{split}
	\Pb\left( \norm{\ED{u}(0,\cdot) - \ED{\phi}}_{L^2(K \times \B^d)} \leq \tilde{\epsilon} \right)
	&\geq 1 - \delta
\end{split}
\end{equation*}%
provided $M$ and $N$ are so large that $N \geq 16cT$ and
\begin{equation*}
	\frac{\tilde{\epsilon} \sqrt{\delta}}{M^{-1/2} + \sqrt{\delta} N^{-\alpha}}
	\geq d c_{r,\alpha} \kappa (1+c) T^3 e^{ \kappa c T} \norm{1 + \norm{x}^{5r+10}}_{L^2(K \times \B^d)}.
\end{equation*}%
Now, by \cref{covert_derivative_error},
\begin{equation*}
\begin{split}
	\norm{\nabla u(0,\cdot) - \nabla\phi}_{L^2(K)}^2
	&\leq \frac{2 d}{\pi^{d/2}} \Gamma\left( \frac{d+4}{2} \right) \int_K \int_{\B^d} \norm{\nabla u(0,x_1)^Tx_2 - \nabla\phi(x_1)^Tx_2}^2 dx_2dx_1
\end{split}
\end{equation*}%
and, trivially,
\begin{equation*}
\begin{split}
	\norm{u(0,\cdot) - \phi}_{L^2(K)}^2
	&= \frac{1}{\mathrm{Vol}(\B^d)} \int_K \int_{\B^d} \norm{u(0,x_1) - \phi(x_1)}^2 dx_2 dx_1
	\\
	&\leq \frac{2 d}{\pi^{d/2}} \Gamma\left( \frac{d+4}{2} \right) \int_K \int_{\B^d} \norm{u(0,x_1) - \phi(x_1)}^2 dx_2 dx_1.
\end{split}
\end{equation*}%
Hence,
\begin{equation*}
\begin{split}
	\norm{u(0,\cdot) - \phi}_{L^2(K)}^2 + \norm{\nabla u(0,\cdot) - \nabla\phi}_{L^2(K)}^2
	&\leq \frac{2 d}{\pi^{d/2}} \Gamma\left( \frac{d+4}{2} \right) \int_{K \times \B^d} \norm{\ED{u}(0,x) - \ED{\phi}(x)}^2 dx.
\end{split}
\end{equation*}%
Then,
\begin{equation*}
\begin{split}
	&\Pb\left( \norm{u(0,\cdot) - \phi}_{L^2(K)} + \norm{\nabla u(0,\cdot) - \nabla\phi}_{L^2(K)} \leq \epsilon \right)
	\\
	&\geq \Pb\left( \norm{\ED{u}(0,\cdot) - \ED{\phi}}_{L^2(K \times \B^d)} \leq \frac{\epsilon}{\sqrt{\frac{4 d}{\pi^{d/2}} \Gamma\left( \frac{d+4}{2} \right)}} \right) \geq 1-\delta
\end{split}
\end{equation*}%
provided $M$ and $N$ are so large that
\begin{equation*}
	\frac{\epsilon \sqrt{\delta}}{M^{-1/2} + \sqrt{\delta} N^{-\alpha}}
	\geq d c_{r,\alpha} \kappa (1+c) T^3 e^{ \kappa c T} \norm{1 + \norm{x}^{5r+10}}_{L^2(K \times \B^d)} \sqrt{\frac{4 d}{\pi^{d/2}} \Gamma\left( \frac{d+4}{2} \right)}.
\end{equation*}%
Abbreviating $s = (5r+10)/2$ and using that $(1+(\alpha^2+\beta^2)^s)^2 \leq 2^{2s-1}(1+\alpha^{2s})^2(1+\beta^{4s})$, we find
\begin{equation*}
\begin{split}
	\norm{1 + \norm{x}^{5r+10}}_{L^2(K \times \B^d)}^2
	&= \int_K \int_{\B^d} \left( 1 + \left( \norm{x_1}^2 + \norm{x_2}^2 \right)^s \right)^2 dx_2dx_1
	\\
	&\leq 2^{2s-1} \norm{1 + \norm{x}^{5r+10}}_{L^2(K)}^2 \int_{\B^d} \left( 1 + \norm{x_2}^{4s} \right) dx_2.
\end{split}
\end{equation*}%
Moreover,
\begin{equation*}
\begin{split}
	\int_{\B^d} \left( 1 + \norm{x_2}^{4s} \right) dx_2
	&= \mathrm{Vol}(\partial\B^d) \int_0^1 \left( 1 + \tau^{4s} \right) \tau^{d-1} d\tau
	\\
	&= \mathrm{Vol}(\partial\B^d) \frac{4s+2d}{(4s+d)d}
	= \frac{2\pi^{d/2}}{\Gamma\left( \frac{d}{2} \right)} \frac{4s+2d}{(4s+d)d}.
\end{split}
\end{equation*}%
Thus,
\begin{equation*}
\begin{split}
	\kappa \norm{1 + \norm{x}^{5r+10}}_{L^2(K \times \B^d)} \sqrt{\frac{4 d}{\pi^{d/2}} \Gamma\left( \frac{d+4}{2} \right)}
	&\leq \kappa \norm{1 + \norm{x}^{5r+10}}_{L^2(K)} \sqrt{ \frac{4s+2d}{4s+d} \frac{\Gamma\left( \frac{d+4}{2} \right)}{\Gamma\left( \frac{d}{2} \right)}}
	\\
	&\leq d \kappa \norm{1 + \norm{x}^{5r+10}}_{L^2(K)},
\end{split}
\end{equation*}%
where we used that $4 \Gamma\left( \frac{d+4}{2} \right) = d(d+2) \Gamma\left( \frac{d}{2} \right)$.
We conclude that a sufficient bound for $M$ and $N$ is
\begin{equation*}
	\frac{\epsilon \sqrt{\delta}}{M^{-1/2} + \sqrt{\delta} N^{-\alpha}}
	\geq d^2 c_{r,\alpha} \kappa (1+c) T^3 e^{ \kappa c T} \norm{1 + \norm{x}^{5r+10}}_{L^2(K)}.
\end{equation*}%
\end{proof}

\begin{remark}
\label{MCES_Sobolev_time_dependent}
	As in \cref{MCES_time_dependent}, if we consider the Euler scheme starting at time $t \in [0,T]$, then we obtain a time-dependent approximation $\phi(t,x)$ with
\begin{equation*}
\begin{split}
	\Pb\left( \norm{u - \phi}_{L^2([0,T] \times K)} + \norm{\nabla u - \nabla\phi}_{L^2([0,T] \times K)} \leq \epsilon \right)
	&\geq 1 - \delta
\end{split}
\end{equation*}%
provided $M$ and $N$ are so large that
\begin{equation*}
	\frac{\epsilon \sqrt{\delta}}{M^{-1/2} + \sqrt{\delta} N^{-\alpha}}
	\geq d^2 c_{r,\alpha} \kappa_r (1+c) T^4 e^{ \kappa_r c T} \norm{1 + \norm{x}^{5r+10}}_{L^2(K)}.
\end{equation*}%
\end{remark}

As in \cref{MCES_error}, the curse of dimensionality is broken if $c_{r,\alpha}$ grows at most polynomially in the dimension and $c$ at most logarithmically.
This time, these constants measure the spatial growth of the augmented derivatives of the coefficients.
To give an example of how $\GrowthSingle{\ED{\mu}}$ might behave, we revisit \cref{ex_PG}.

\begin{example}
\label{ex_PG_comDer}
	Let $\fct \colon \R^d \rightarrow \R^{\dout}$ be differentiable and Lipschitz continuous with Lipschitz constant $L \geq 0$.
Then, for all $r \in [1,\infty)$,
\begin{equation*}
	\GrowthSum{\ED{\fct}}{r}{0}
	\leq \max\left\{ \norm{\fct(0)} , \frac{L}{r} \right\}.
\end{equation*}%
\end{example}

Thus, if the Lipschitz constant of $\mu$ and $\sigma$ grows at most logarithmically in the dimension, then the constant $e^{\kappa_r c T}$ in \cref{Sobolev_MCES_error} grows at most polynomially in the dimension.
Some more specific examples are as follows.

\begin{example}
	Let $\fct_i \colon \R \rightarrow \R$, $0 \leq i \leq d$, be differentiable and Lipschitz continuous with Lipschitz constant $L \geq 0$.
\begin{enumerate}[\rm (i)]\itemsep = 0em
\item
Consider the function $\fct \colon \R^d \rightarrow \R$ given by $\fct(x) = \fct_0(\frac{1}{d}\sum_{i=1}^d \fct_i(x_i))$.
Then, for all $r \in [1,\infty)$,
\begin{equation*}
	\GrowthSum{\ED{\fct}}{r}{0}
	\leq \max\left\{ (L+1) \max_{0 \leq i \leq d} \abs{\fct_i(0)} , \frac{L^2}{r \sqrt{d}} \right\}.
\end{equation*}%
\item
Consider the function $\fct \colon \R^d \rightarrow \R$ given by $\fct(x) = \fct_0(\frac{1}{\sqrt{d}}\sum_{i=1}^d \fct_i(x_i))$ and suppose $\fct_i(0) = 0$ for all $1 \leq i \leq d$.
Then, for all $r \in [1,\infty)$,
\begin{equation*}
	\GrowthSum{\ED{\fct}}{r}{0}
	\leq \max\left\{ \abs{\fct_0(0)} , \frac{L^2}{r} \right\}.
\end{equation*}%
\end{enumerate}	
\end{example}


\section{Monte Carlo Euler scheme with perturbed coefficients}
\label{section_perturbed}

\subsection{Error rate for the perturbed scheme}

In the previous two sections, we derived (Sobolev) error rates for the Monte Carlo Euler scheme.
The next step is to perturb the coefficients $\mu$ and $\sigma$ and the functions $f$ and $g$ and see how the resulting scheme performs.
First, we analyze growth and Lipschitz properties of the Euler scheme, jointly in the initial condition and the coefficients.
In \cref{Euler_scheme_well_defined}, we considered the norm of the Euler scheme as a stochastic process.
In the lemma below, we consider its Euclidean norm point-wise.
Throughout this section, $F \colon \R^d \rightarrow \R^{\dout}$ and $G \colon [0,T] \times \R^d \rightarrow \R^{\dout}$ denote functions and $\nu \colon [0,T] \times \R^d \rightarrow \R^d$ and $\tau \colon [0,T] \times \R^d \rightarrow \R^{d \times \dWiener}$ denote measurable functions, which we think of as perturbations of $f$, $g$, $\mu$, and $\sigma$.

\begin{lemma}
\label{Euler_perturbation}
	Denote $c_0 = \max\{ \GrowthSingle{\mu} , \GrowthSingle{\sigma} \}$ and $\Delta = \frac{T}{N} + \max_{0 \leq n \leq N-1} \norm{ W_{t_{n+1}}-W_{t_n} }$.
Then, the following hold.
\begin{enumerate}[\rm (i)]\itemsep = 0em
\item
We have, almost surely for all $0 \leq n \leq N$ and $x \in \R^d$,
\begin{equation}
\label{Euler_perturbation_item_one}
\begin{split}
	1 + \norm{ \Euler(x)_{t_n} }
	&\leq (1+ c_0 \Delta)^n (1+\norm{x}).
\end{split}
\end{equation}%
\item
If, in addition, $\mu$ and $\sigma$ are differentiable in the space variable, then, almost surely for all $0 \leq n \leq N$ and $x,y \in \R^d$,
\begin{equation}
\label{Euler_perturbation_item_two}
\begin{split}
	\norm{ \Euler(x)_{t_n} - \EulerPert(y)_{t_n} }
	&\leq (1 + c_1 \Delta)^{nN} (1 + \max\{\norm{x},\norm{y}\})^n \left( \norm{x-y} + \eta \right),
\end{split}
\end{equation}%
where $c_1 = \max\{ 1 , \GrowthSingle{\mu} , \GrowthSingle{\nabla\mu} , \GrowthSingle{\sigma} , \GrowthSingle{\nabla\sigma} , \GrowthSingle{\nu} , \GrowthSingle{\tau} \}$ and $\eta = \max\{ \GrowthSingle{\mu-\nu} , \GrowthSingle{\sigma-\tau} \}$.
\end{enumerate}
\end{lemma}

\begin{proof}
	Fix $x,y \in \R^d$ and abbreviate $\Ec_n = \Euler(x)_{t_n}$ and $\Ec'_n = \EulerPert(y)_{t_n}$.
We first show \eqref{Euler_perturbation_item_one}.
To this end, observe that $\norm{ \mu(t_n,\Ec_n) } \leq \GrowthSingle{\mu} (1+\norm{\Ec_n})$ and likewise for $\sigma$.
The recursion
\begin{equation*}
\begin{split}
	\Ec_{n+1}
	&= \Ec_n + \mu(t_n,\Ec_n) (t_{n+1}-t_n) + \sigma(t_n,\Ec_n) (W_{t_{n+1}}-W_{t_n})
\end{split}
\end{equation*}%
translates into the bound
\begin{equation*}
\begin{split}
	1 + \norm{ \Ec_{n+1} }
	&\leq 1 + \norm{ \Ec_n } + c_0 (1 + \norm{ \Ec_n }) \Delta
	= (1+c_0\Delta) (1+\norm{ \Ec_n }).
\end{split}
\end{equation*}%
With $\norm{ \Ec_0 } = \norm{x}$, we immediately deduce \eqref{Euler_perturbation_item_one}.
The inequality \eqref{Euler_perturbation_item_two} is shown similarly once we note that
\begin{equation*}
\begin{split}
	\norm{ \mu(t_n,\Ec_n) - \mu(t_n,\Ec'_n) }
	&\leq \GrowthSingle{\nabla\mu} (1+\max\{\norm{\Ec_n},\norm{\Ec'_n}\}) \norm{\Ec_n-\Ec'_n}
\end{split}
\end{equation*}%
and
\begin{equation*}
\begin{split}
	\norm{ \mu(t_n,\Ec'_n) - \nu(t_n,\Ec'_n) }
	&\leq \GrowthSingle{\mu-\nu} (1+\norm{\Ec'_n}).
\end{split}
\end{equation*}%
Analogous bounds hold for $\sigma$ and $\tau$.
Plugging this into the recursion defining $\Ec_n$ and $\Ec'_n$ and abbreviating $C = \max\left\{ \GrowthSingle{\nabla\mu} , \GrowthSingle{\nabla\sigma} \right\}$ and $R = 1 + \max_{0 \leq n \leq N-1} \max\{\norm{\Ec_n},\norm{\Ec'_n}\}$, we find
\begin{equation*}
\begin{split}
	\norm{ \Ec_{n+1} - \Ec'_{n+1} }
	&\leq \norm{ \Ec_n - \Ec'_n} + \norm{ \mu(t_n,\Ec_n) - \nu(t_n,\Ec'_n)} (t_{n+1}-t_n) + \norm{ \sigma(t_n,\Ec_n) - \tau(t_n,\Ec'_n) } \norm{ W_{t_{n+1}}-W_{t_n} }
	\\
	&\leq \norm{ \Ec_n - \Ec'_n} (1 + C R \Delta) + \eta R \Delta.
\end{split}
\end{equation*}%
Since $\norm{\Ec_0-\Ec'_0}=\norm{x-y}$, expanding this recursion yields
\begin{equation*}
\begin{split}
	\norm{\Ec_n-\Ec'_n}
	&\leq \norm{x-y} (1 + C R \Delta)^n + \eta R \Delta \sum_{k=0}^{n-1} (1 + C R \Delta)^k
	\\
	&= \norm{x-y} (1 + C R \Delta)^n + \eta \frac{(1 + C R \Delta)^n - 1}{C}
	\leq \left( \norm{x-y} + \eta \right) (1 + c_1 R \Delta)^n.
\end{split}
\end{equation*}%
We proved in \eqref{Euler_perturbation_item_one} that $R \leq (1 + c_1\Delta)^{N-1} (1 + \max\{\norm{x},\norm{y}\})$, with which we conclude
\begin{equation*}
\begin{split}
	1 + c_1 R \Delta
	&\leq (1 + c_1 \Delta)^N (1 + \max\{\norm{x},\norm{y}\}).
\end{split}
\end{equation*}%
\end{proof}

We translate the bound from the previous lemma into the expectation appearing in the Feynman--Kac formula.
Then, we will apply Monte Carlo sampling to the expectation in the function $w$ specified below instead of to the expectation in $v$ as done previously.

\begin{proposition}
\label{Euler_Exp_perturbation}
	Assume $\mu$, $\sigma$, $f$, $g$, $F$, and $G$ are differentiable in the space variable.
Consider $v \colon \R^d \rightarrow \R^{\dout}$ and $w \colon \R^d \rightarrow \R^{\dout}$ given by
\begin{equation*}
\begin{split}
	v(x)
	&= \Exp{ f\left(\Euler(x)_T\right) + \frac{T}{N} \sum_{n=0}^{N-1} g\left(t_n,\Euler(x)_{t_n}\right) },
	\\
	w(x)
	&= \Exp{ F\left(\EulerPert(x)_T\right) + \frac{T}{N} \sum_{n=0}^{N-1} G\left(t_n,\EulerPert(x)_{t_n}\right) }.
\end{split}
\end{equation*}%
Denote
\begin{equation*}
\begin{split}
	c_0
	&= \max\left\{ 1 , \GrowthSingle{\nabla f} , \GrowthSingle{\nabla g} \right\},
	\\
	c_1
	&= \max\left\{ 1 , \GrowthSingle{\mu} , \GrowthSingle{\nabla\mu} , \GrowthSingle{\sigma} , \GrowthSingle{\nabla\sigma} , \GrowthSingle{\nu} , \GrowthSingle{\tau} \right\},
	\\
	\eta
	&= \max\left\{ \GrowthSingle{\mu-\nu} , \GrowthSingle{\sigma-\tau} , \GrowthSingle{f-F} , \GrowthSingle{g-G} \right\}.
\end{split}
\end{equation*}%
If $N \geq \dWiener+T+1$, then, for all $x \in \R^d$,
\begin{equation*}
\begin{split}
	\norm{v(x)-w(x)}
	&\leq \frac{3}{4} \eta c_0 T N \left( 6 c_1 \sqrt{T N} \right)^{N(1+N)} \left(1 + \norm{x}^{N+1}\right).
\end{split}
\end{equation*}%
\end{proposition}

\begin{proof}
	Abbreviate $\Ec = \Euler(x)$ and $\Ec' = \EulerPert(x)$.
As in the proof of \cref{Euler_perturbation},
\begin{equation*}
\begin{split}
	\norm{ f(\Ec(x)_T) - F(\Ec'(x)_T) }
	&\leq \norm{ f(\Ec(x)_T) - f(\Ec'(x)_T) } + \norm{ f(\Ec'(x)_T) - F(\Ec'_x)_T) }
	\\
	&\leq \GrowthSingle{\nabla f} (1+\max\{\norm{\Ec(x)_T},\norm{\Ec'(x)_T}\}) \norm{\Ec(x)_T-\Ec'(x)_T}
	\\
	&\quad + \GrowthSingle{f-F} (1+\norm{\Ec'(x)_T}).
\end{split}
\end{equation*}%
Abbreviate $\Delta = \frac{T}{N} + \max_{0 \leq n \leq N-1} \norm{ W_{t_{n+1}}-W_{t_n} }$.
Applying \cref{Euler_perturbation}, we obtain
\begin{equation*}
\begin{split}
	\norm{ f(\Ec(x)_T) - F(\Ec'(x)_T) }
	&\leq \eta c_0 (1+ c_1 \Delta)^{N(1+N)} (1 + \norm{x})^{N+1} + \eta (1+ c_1 \Delta)^N (1+\norm{x})
	\\
	&\leq 2 \eta c_0 (1+ c_1 \Delta)^{N(1+N)} (1 + \norm{x})^{N+1}.
\end{split}
\end{equation*}%
The same bound holds for $\norm{ g(t_n,\Ec(x)_{t_n}) - G(t_n,\Ec'(x)_{t_n}) }$.
Thus,
\begin{equation*}
\begin{split}
	\norm{v(x)-w(x)}
	&\leq \Exp{ \norm{f(\Ec(x)_T)-F(\Ec'(x)_T)} + \frac{T}{N} \sum_{n=0}^{N-1} \norm{g(t_n,\Ec(x)_{t_n})-G(t_n,\Ec'(x)_{t_n})} }
	\\
	&\leq 4 \eta c_0 T \Exp{ (1+ c_1 \Delta)^{N(1+N)} } (1 + \norm{x})^{N+1}.
\end{split}
\end{equation*}%
Moreover, with $k = N(1+N)$,
\begin{equation*}
\begin{split}
	\Exp{ (1+ c_1 \Delta)^k }
	&\leq 3^{k-1} c_1^k \Exp{ 1 + \left(\frac{T}{N}\right)^k + \left( \Delta - \frac{T}{N} \right)^k }
	\leq 3^{k-1} c_1^k \Exp{ 2 + \left( \Delta - \frac{T}{N} \right)^k },
\end{split}
\end{equation*}%
where we used that $T \leq N$.
If $Z \sim \Nc(0,(T/N) \cdot \Ind_{\dWiener})$ is multivariate normal, then
\begin{equation*}
\begin{split}
	\Exp{ \left( \Delta - \frac{T}{N} \right)^k }
	\leq N \Exp{ \norm{Z}^k }
	&= N \left(\frac{2 T}{N}\right)^{k/2} \frac{\Gamma\left(\frac{\dWiener+k}{2}\right)}{\Gamma\left(\frac{\dWiener}{2}\right)}
	\\
	&\leq N \left(\frac{T (\dWiener+k-2)}{N}\right)^{k/2}
	\leq N \left(\frac{2 T k}{N}\right)^{k/2},
\end{split}
\end{equation*}%
since $\dWiener \leq N \leq k$.
Thus,
\begin{equation*}
\begin{split}
	\Exp{ (1+ c_1 \Delta)^k }
	&\leq 3^{k-1} c_1^k \left( 2 + N \left(\frac{2 T k}{N}\right)^{k/2} \right)
	\leq \frac{2}{3} N c_1^k \left(\frac{18 T k}{N}\right)^{k/2}.
\end{split}
\end{equation*}%
In total,
\begin{equation*}
\begin{split}
	\norm{v(x)-w(x)}
	&\leq \frac{8}{3} \eta c_0 T N \left( c_1 \sqrt{24 T N} \right)^{N(1+N)} (1 + \norm{x})^{N+1}
	\\
	&\leq \frac{3}{4} \eta c_0 T N \left( 6 c_1 \sqrt{T N} \right)^{N(1+N)} \left(1 + \norm{x}^{N+1}\right).
\end{split}
\end{equation*}%
\end{proof}

Now, we immediately obtain the error estimate for the perturbed Monte Carlo Euler scheme.

\begin{corollary}
\label{MCES_perturbation}
	Let \cref{assumption_weak_error} hold and suppose $\nu$, $\tau$, $F$, and $G$ are differentiable in the space variable.
Let $\Ec_1(x),\dots,\Ec_M(x)$ be $M \in \N$ independent copies of $\EulerPert(x)$.
Consider $\psi \colon \Omega \times \R^d \rightarrow \R^{\dout}$ given by
\begin{equation*}
\begin{split}
	\psi(x) &= \frac{1}{M} \sum_{m=1}^M \left[ F\left(\Ec_m(x)_T\right) + \frac{T}{N} \sum_{n=0}^{N-1} G\left(t_n,\Ec_m(x)_{t_n}\right) \right].
\end{split}
\end{equation*}%
Let $\alpha \in (0,1]$, $\delta,\epsilon \in (0,1)$, $r \in [0,\infty)$, let $K \subseteq \R^d$ be compact, and denote
\begin{equation*}
\begin{split}
	c
	&= \max\{\GrowthSingle{\mu},\GrowthSingle{\sigma}^2,\GrowthSingle{\nu},\GrowthSingle{\tau}^2\},
	\\
	c_{r,\alpha}
	&= \max\big\{ \GrowthSum{F}{r}{0} , \GrowthSum{G}{r}{0} , \Cc_{r,\alpha}(\mu,\sigma,g,u) \big\},
	\\
	c_0
	&= \max\left\{ 1 , \GrowthSingle{\nabla f} , \GrowthSingle{\nabla g} \right\},
	\\
	c_1
	&= \max\left\{ 1 , \GrowthSingle{\mu} , \GrowthSingle{\nabla\mu} , \GrowthSingle{\sigma} , \GrowthSingle{\nabla\sigma} , \GrowthSingle{\nu} , \GrowthSingle{\tau} \right\},
	\\
	\eta
	&= \max\left\{ \GrowthSingle{\mu-\nu} , \GrowthSingle{\sigma-\tau} , \GrowthSingle{f-F} , \GrowthSingle{g-G} \right\};
\end{split}
\end{equation*}%
see \eqref{cumbersome_constant}.
If $N \geq 16cT$ and $N \geq \dWiener+T+1$, then there exists a constant $\kappa_r > 0$, depending only on $r$, such that
\begin{equation*}
\begin{split}
	\Pb\left( \norm{u(0,\cdot) - \psi}_{L^2(K)} \leq \epsilon \right)
	&\geq 1 - \delta
\end{split}
\end{equation*}%
provided $M$ and $N$ are so large that
\begin{equation*}
	\frac{\epsilon \sqrt{\delta}}{M^{-1/2} + \sqrt{\delta} N^{-\alpha}}
	\geq d c_{r,\alpha} \kappa_r (1+c) T^3 e^{ \kappa_r c T} \norm{1 + \norm{x}^{5r+10}}_{L^2(K)}
\end{equation*}%
and $\eta$ is so small that
\begin{equation*}
	\frac{\epsilon}{\eta}
	\geq c_0 T N \left( 6 c_1 \sqrt{T N} \right)^{N(1+N)} \norm{ 1 + \norm{x}^{N+1} }_{L^2(K)}.
\end{equation*}%
\end{corollary}

\begin{proof}
	This is immediate from \cref{weak_error_rate_thrm,MC_error,Euler_Exp_perturbation} (we apply \cref{MC_error} with $r+1$ to the Euler scheme $\EulerPert$ with coefficients $\nu$ and $\tau$ and to the functions $F$ and $G$).
\end{proof}

We see in \cref{MCES_perturbation} that the requirement on $M$ and $N$ is verbatim from \cref{MCES_error}.
In particular, this does not incur the curse of dimensionality if we can control $c$ and $c_{r,\alpha}$.
However, this is not true for the requirement on $\eta$, which scales exponentially in $N$.
The error inferred from perturbing the coefficients is harder to control than the error from the Monte Carlo Euler scheme itself.
As such, the perturbed scheme is provably useful only if the cost of the perturbations in a numerical sense is exponentially small.


\subsection{Sobolev error rate for the perturbed scheme}

As with \cref{Sobolev_MCES_error}, the Sobolev error rate for the perturbed Monte Carlo Euler scheme is obtained by applying \cref{MCES_perturbation} to the augmented derivatives and using the trick from \cref{covert_derivative_error}.
In particular, we need to impose one more order of regularity on all functions involved.

\begin{assumption}
\label{Sobolev_perturbation_assumption}
	Let \cref{Sobolev_assumption_weak_error} hold and suppose $\nu$, $\tau$, $F$, and $G$ are twice differentiable in the space variable.
Suppose that $\nu$ and $\tau$ also satisfy \cref{assumption_SDE_diff}.
\end{assumption}

\begin{theorem}
\label{MCES_Sobolev_perturbation}
	Let \cref{Sobolev_perturbation_assumption} hold and let $\Ec_1(x),\dots,\Ec_M(x)$ be $M \in \N$ independent copies of $\EulerPert(x)$.
Consider $\psi \colon \Omega \times \R^d \rightarrow \R^{\dout}$ given by
\begin{equation*}
\begin{split}
	\psi(x) &= \frac{1}{M} \sum_{m=1}^M \left[ F\left(\Ec_m(x)_T\right) + \frac{T}{N} \sum_{n=0}^{N-1} G\left(t_n,\Ec_m(x)_{t_n}\right) \right].
\end{split}
\end{equation*}%
Let $\alpha \in (0,1]$, $\delta,\epsilon \in (0,1)$, $r \in [0,\infty)$, let $K \subseteq \R^d$ be compact, and denote
\begin{equation*}
\begin{split}
	c
	&= \max\{\GrowthSingle{\ED{\mu}},\GrowthSingle{\ED{\sigma}}^2,\GrowthSingle{\ED{\nu}},\GrowthSingle{\ED{\tau}}^2\},
	\\
	c_{r,\alpha}
	&= \max\{ \GrowthSum{\ED{F}}{r}{0} , \GrowthSum{\ED{G}}{r}{0} , \Cc_{r,\alpha}(\ED{\mu},\ED{\sigma},\ED{g},\ED{u}) \},
	\\
	c_0
	&= \max\{ 1 , \GrowthSingle{\nabla \ED{f}} , \GrowthSingle{\nabla \ED{g}} \},
	\\
	c_1
	&= \max\left\{ 1 , \GrowthSingle{\ED{\mu}} , \GrowthSingle{\nabla\ED{\mu}} , \GrowthSingle{\ED{\sigma}} , \GrowthSingle{\nabla\ED{\sigma}} , \GrowthSingle{\ED{\nu}} , \GrowthSingle{\ED{\tau}} \right\},
	\\
	\eta
	&= \max\left\{ \GrowthSingle{\WED{\mu-\nu}} , \GrowthSingle{\WED{\sigma-\tau}} , \GrowthSingle{\WED{f-F}} , \GrowthSingle{\WED{g-G}} \right\};
\end{split}
\end{equation*}%
see \eqref{cumbersome_constant}.
If $N \geq 16cT$ and $N \geq \dWiener+T+1$, then there exists a constant $\kappa_r > 0$, depending only on $r$, such that
\begin{equation*}
\begin{split}
	\Pb\left( \norm{u(0,\cdot) - \psi}_{L^2(K)} + \norm{\nabla u(\cdot,0) - \nabla\psi}_{L^2(K)} \leq \epsilon \right)
	&\geq 1 - \delta
\end{split}
\end{equation*}%
provided $M$ and $N$ are so large that
\begin{equation*}
	\frac{\epsilon \sqrt{\delta}}{M^{-1/2} + \sqrt{\delta} N^{-\alpha}}
	\geq d^2 c_{r,\alpha} \kappa_r (1+c) T^3 e^{ \kappa_r c T} \norm{1 + \norm{x}^{5r+10}}_{L^2(K)}
\end{equation*}%
and $\eta$ is so small that
\begin{equation*}
	\frac{\epsilon}{\eta}
	\geq 3 d c_0 T N \left( 7 c_1 \sqrt{T N} \right)^{N(1+N)} \norm{ 1 + \norm{x}^{N+1} }_{L^2(K)}.
\end{equation*}%
\end{theorem}

\begin{proof}
	We deduce this from \cref{MCES_perturbation} the same way we deduced \cref{Sobolev_MCES_error} from \cref{MCES_error} with the additional analogous computation
\begin{equation*}
\begin{split}
	\norm{ 1 + \norm{x}^{N+1} }_{L^2(K \times \B^d)} \sqrt{\frac{4 d}{\pi^{d/2}} \Gamma\left( \frac{d+4}{2} \right)}
	&\leq 2^{(N+5)/2} d \norm{ 1 + \norm{x}^{N+1} }_{L^2(K)}.
\end{split}
\end{equation*}%
The requirement on $\nu$ and $\tau$ to satisfy \cref{assumption_SDE_diff} ensures that $\EulerPert$ commutes with the augmented derivative by \cref{comDer_Euler_commute}.
\end{proof}


\subsection{Neural network approximations}

\begin{figure}
\begin{center}
\includegraphics[scale=0.85]{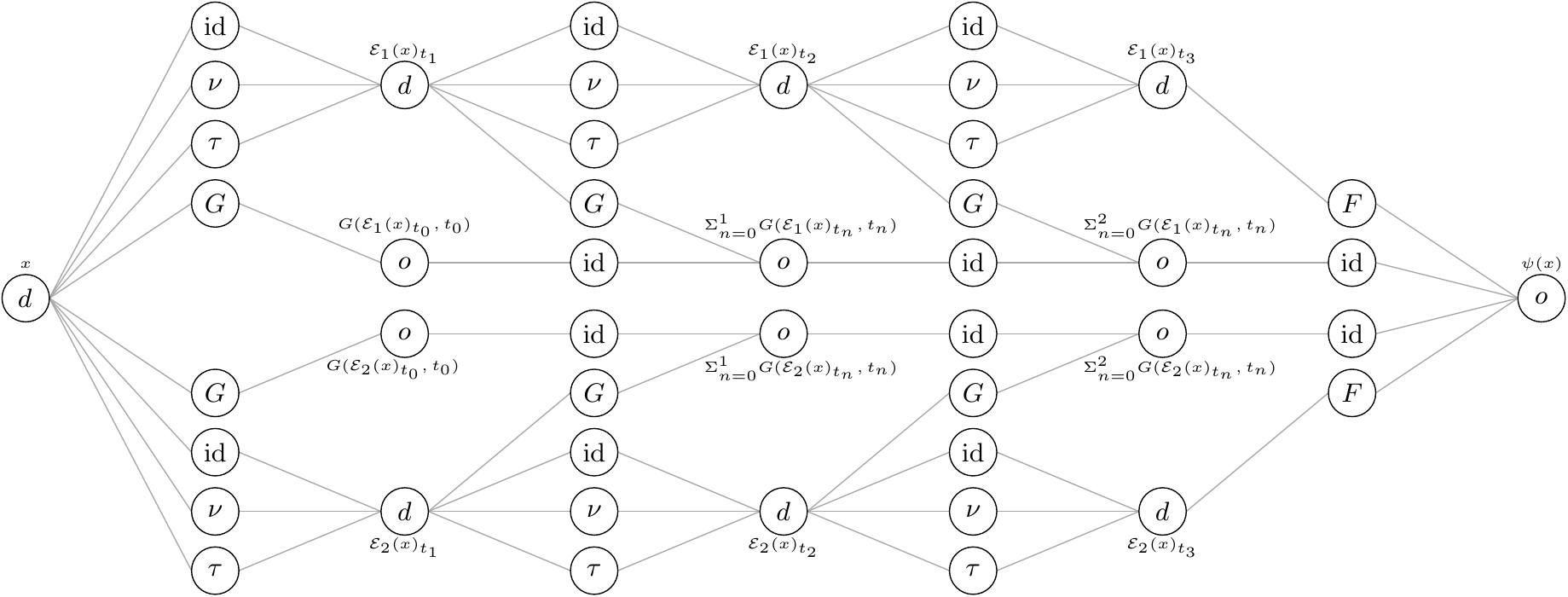}
\caption{Network implementing the Monte Carlo Euler scheme in \cref{MCES_NN_implementation} with $M=2$ and $N=3$.
Each circle represents itself a small network.}
\label{fig_MCES_NN}
\end{center}
\end{figure}

We mentioned in the introduction that neural networks are particularly handy since their compositional structure makes them well-suited to imitate the Monte Carlo Euler scheme.
A network is a concatenation $\Ac_D \circ \rho \circ \dots \circ \Ac_2 \circ \rho \circ \Ac_1$ of affine functions $\Ac_j$ and a fixed nonlinearity $\rho \colon \R \rightarrow \R$, which is applied component-wise.
We say the network has depth $D \in \N$ and consists of $D+1$ layers.
The tuple of the input, respectively output dimensions of the affine functions comprise the architecture of the network.
We allow for identity nodes, meaning that one may also apply the identity function instead of the nonlinearity $\rho$ at any point in the network.
The parameters that characterize the affine functions are called the parameters of the network.
Let $\Pc(\cdot)$ denote the number of non-zero parameters.
The next lemma showcases that networks can implement the Monte Carlo Euler scheme efficiently if the coefficients are given by networks.
It is straight-forward to prove; see \cref{fig_MCES_NN} for an illustration of the architecture.

\begin{assumption}
\label{assumption_NNs}
Suppose, for each fixed $t \in [0,T]$, that the functions $F$, $G(t,\cdot)$, $\nu(t,\cdot)$, and $\tau(t,\cdot)$ are given by networks, whose architectures are independent of $t$ and whose depths are all the same, namely $D \in \N$.
\end{assumption}

\begin{lemma}
\label{MCES_NN_implementation}
	Let \cref{assumption_NNs} hold and let $\Ec_1(x),\dots,\Ec_M(x)$ be $M \in \N$ copies of $\EulerPert(x)$.
Then, almost surely for all $M,N \in \N$, the function
\begin{equation*}
\begin{split}
	x \mapsto \psi(x)
	&= \frac{1}{M} \sum_{m=1}^M \left[ F\left(\Ec_m(x)_T\right) + \frac{T}{N} \sum_{n=0}^{N-1} G\left(t_n,\Ec_m(x)_{t_n}\right) \right]
\end{split}
\end{equation*}%
can be given by a network with
\begin{equation*}
\begin{split}
	\Pc(\psi)
	&\leq M ( \Pc(F) + N [ \Pc(G) + \Pc(\nu) + \Pc(\tau) + (d+\dout)D ] ).
\end{split}
\end{equation*}%
\end{lemma}

Finally, we use networks for Sobolev approximations of the PDE solution $u$.
To this end, we assume that the coefficients can be efficiently approximated with networks and then use these networks in the perturbed scheme.
\cref{MCES_Sobolev_perturbation} asserts that there exists an element $\omega \in \Omega$ such that a realization of the perturbed scheme approximates $u$ to accuracy $\epsilon$.
This realization can be implemented by a network as shown in the previous lemma.
After \cref{MCES_perturbation}, we remarked that the cost of the perturbations needs to be exponentially small.
This means that the coefficients need to be approximated by the respective networks at an algebraic rate.

\begin{theorem}
\label{NN_error}
	Let \cref{Sobolev_perturbation_assumption,assumption_NNs} hold.
Let $\alpha \in (0,1]$, $\epsilon \in (0,1)$, $r \in [0,\infty)$, and denote
\begin{equation*}
\begin{split}
	c
	&= \max\{\GrowthSingle{\ED{\mu}},\GrowthSingle{\ED{\sigma}}^2,\GrowthSingle{\ED{\nu}},\GrowthSingle{\ED{\tau}}^2\}.
\end{split}
\end{equation*}%
Suppose that the quantities
\begin{equation*}
\begin{split}
	&c , \GrowthSum{\mu}{0}{2} , \GrowthSum{\sigma}{0}{2} , \GrowthSum{f}{0}{2} , \GrowthSum{g}{0}{2} ,
	\\
	&\GrowthSum{F}{r}{1} , \GrowthSum{G}{r}{1} , \GrowthSum{\partial_t g}{r}{1} , \GrowthSum{\mu}{r}{3} , \GrowthSum{\sigma}{r}{3} , \GrowthSum{g}{r}{3} , \GrowthSum{u}{r}{5} , \GrowthHoelder{\ED{\mu}}{r}{\alpha} , \GrowthHoelder{\ED{\sigma}}{r}{\alpha}
\end{split}
\end{equation*}%
are all bounded by $\kappa_{r,\alpha} d^{\kappa_{r,\alpha}}$ for a constant $\kappa_{r,\alpha}$ depending only on $r$ and $\alpha$, and suppose
\begin{equation*}
\begin{split}
	\max\left\{ \GrowthSum{\mu-\nu}{0}{1} , \GrowthSum{\sigma-\tau}{0}{1} , \GrowthSum{f-F}{0}{1} , \GrowthSum{g-G}{0}{1} \right\}
	&\leq \lambda^{d^{\lambda}} \exp\left( - \Pc^{\lambda^{-1}} \right),
\end{split}
\end{equation*}%
where $\lambda$ is a universal constant and $\Pc = \min\{ \Pc(F) , \Pc(G) , \Pc(\nu) , \Pc(\tau) \}$.
Then, there exists a constant $\kappa > 0$, depending only on $T$, $\lambda$, $r$, and $\alpha$, and a network $\psi$ such that
\begin{equation*}
\begin{split}
	\norm{u(0,\cdot) - \psi}_{L^2([0,1]^d)} + \norm{\nabla u(0,\cdot) - \nabla\psi}_{L^2([0,1]^d)} \leq \epsilon
\end{split}
\end{equation*}%
and
\begin{equation*}
\begin{split}
	\Pc(\psi)
	\leq \kappa (d \dout)^{\kappa} e^{\kappa c} \epsilon^{-\kappa}.
\end{split}
\end{equation*}%
\end{theorem}

\begin{proof}
	This follows from
\cref{MCES_Sobolev_perturbation,MCES_NN_implementation}.
That it suffices to require polynomial growth in the dimension of $\GrowthSum{\mu}{0}{2}$ as opposed to $\GrowthSingle{\ED{\mu}}$ and $\GrowthSingle{\nabla\ED{\mu}}$ et cetera follows from the ability to rid the function in $\GrowthSum{\cdot}{s}{k}$ of (augmented) derivatives in exchange for increasing the index $k$ and decreasing the index $s$; see \cref{PGcalc_int} and \cref{PGcalc_Cint} in the appendix.
\end{proof}

We saw in \cref{ex_PG_comDer} that the Lipschitz constant of $\mu$ and $\sigma$ growing at most logarithmically in the dimension is a sufficient criterion to ensure polynomial growth of $e^{\kappa_r c T}$ in the dimension.
For the network approximation, we would require the same from $\nu$ and $\tau$.
Fortunately, network approximations of a Lipschitz continuous function can often be asserted to meet the same Lipschitz continuity as the approximated function;
\cite{CheJenRos2022IEEE}.
In that case, polynomial growth of $e^{\kappa c}$ in the dimension is guaranteed also in \cref{NN_error}.
This concludes the first proved Sobolev approximation of solutions to PDE \eqref{Intro_PDE} with non-affine coefficients by networks free of the curse of dimensionality, up to controlling the constant $c$.


\begin{appendix}

\section{Polynomial growth calculus}

We develop a `polynomial growth calculus' for the growth measure $\GrowthSum{\cdot}{\cdot}{\cdot}$ introduced in \cref{def_PG}.
In particular, we prove \cref{PGcalc_double_generator}.

\subsection{Polynomial growth of derivatives}

The first result shows that if $\GrowthSum{\nabla^k \fct}{r}{0} < \infty$ for some $r \geq 0$, then $\GrowthSum{\fct}{r}{k} < \infty$ with the same $r$.

\begin{lemma}
\label{PGcalc_diff}
	Let $\fct \colon \R^d \rightarrow \R^{\dout}$ be $\N_0 \ni k$-times differentiable and $r \in [0,\infty)$.
Then, for any $0 \leq m \leq k$, we have
\begin{equation*}
	\GrowthSum{\fct}{r}{k}
	\leq \frac{r+k+1}{r+m+1} \GrowthSum{\nabla^{k-m} \fct}{r}{m} + \sum_{j=m+1}^k \frac{r+k+1}{r+j+1} \norm{\nabla^{k-j} \fct(0)}.
\end{equation*}%
\end{lemma}

\begin{proof}
	We claim that
\begin{equation*}
	\GrowthSum{\nabla^{k-j-1} \fct}{r}{j+1}
	\leq \frac{r+j+2}{r+j+1} \GrowthSum{\nabla^{k-j} \fct}{r}{j} + \norm{\nabla^{k-j-1} \fct(0)}
\end{equation*}%
for all $0 \leq j \leq k-1$.
The statement of the lemma is shown by applying this inequality repeatedly.
Now, note that
\begin{equation*}
\begin{split}
	\GrowthSum{\nabla^{k-j-1} \fct}{r}{j+1}
	&\leq \GrowthSum{\nabla^{k-j} \fct}{r}{j} + \sup_{x \in \R^d} \frac{\norm{\nabla^{k-j-1} \fct(x)}}{(1+\norm{x})^{r+j+1}}.
\end{split}
\end{equation*}%
A basic calculus argument shows that
\begin{equation*}
\begin{split}
	\norm{\nabla^{k-j-1} \fct(x) - \nabla^{k-j-1} \fct(0)}
	\leq \int_0^1 \norm{\nabla^{k-j} \fct(sx)} \norm{x} ds
	&\leq \GrowthSum{\nabla^{k-j} \fct}{r}{j} \int_0^1 (1+\norm{sx})^{r+j} \norm{x} ds
	\\
	&\leq \GrowthSum{\nabla^{k-j} \fct}{r}{j} \frac{(1+\norm{x})^{r+j+1}}{r+j+1}.
\end{split}
\end{equation*}%
Thus,
\begin{equation*}
	\sup_{x \in \R^d} \frac{\norm{\nabla^{k-j-1} \fct(x)}}{(1+\norm{x})^{r+j+1}}
	\leq \frac{\GrowthSum{\nabla^{k-j} \fct}{r}{j}}{r+j+1} + \norm{\nabla^{k-j-1} \fct(0)},
\end{equation*}%
from which the claim follows.
\end{proof}

The second result interchanges higher-order derivatives and the augmented derivative.

\begin{lemma}
\label{PGcalc_commute_comDer_stdDer}
	Let $\fct \colon \R^d \rightarrow \R^{\dout}$ be $\N \ni k$-times differentiable and $r \in [0,\infty)$.
Then, for any $0 \leq l,m \leq k$ with $l+m+1 \leq k$, we have
\begin{equation*}
\begin{split}
	\GrowthSum{\WED{\nabla^l \fct}}{r}{m}
	&\leq (m+1) \GrowthSum{\nabla^l \ED{\fct}}{r}{m}
	\\
	\GrowthSum{\nabla^l \ED{\fct}}{r}{m}
	&\leq (m+l+1) \GrowthSum{\WED{\nabla^l \fct}}{r}{m}.
\end{split}
\end{equation*}%
\end{lemma}

\begin{proof}
	We claim for every $0 \leq i \leq k-1$ that $\nabla^i\ED{\fct}(x_1,x_2)$ is a block-tensor that contains $\nabla^i \fct(x_1)$ as exactly $(i+1)$-many blocks, contains $\nabla^{i+1}\fct(x_1)^Tx_2$ as exactly one block, and all of its remaining blocks are zero.
If the claim holds, then, for any $x \in \R^{2d}$,
\begin{equation*}
	\norm{\WED{\nabla^i \fct}(x)}
	\leq \norm{\nabla^i \ED{\fct}(x)}
	\leq (i+1) \norm{\WED{\nabla^i \fct}(x)},
\end{equation*}%
which implies the statement of the lemma.
The claim is immediate from an inductive argument.
The case $i=0$ holds by definition.
For $i \geq 1$, if we differentiate $\nabla^{i-1}\ED{\fct}$ again, we differentiate $i$-many blocks $\nabla^{i-1}\fct(x_1)$ with respect to $x_1$, which yields $i$-many blocks $\nabla^i \fct(x_1)$;
we differentiate the one block $\nabla^i \fct(x_1)^Tx_2$ with respect to $x_2$, which yields the $(i+1)$-th block $\nabla^i \fct(x_1)$;
and we differentiate the one block $\nabla^i \fct(x_1)^Tx_2$ with respect to $x_1$, which yields the one block $\nabla^{i+1}\fct(x_1)^Tx_2$.
\end{proof}

\cref{PGcalc_diff} showed that we can decrease the subscript of $\GrowthSum{\fct}{r}{\cdot}$ by trading derivatives of $\fct$.
Conversely, \cref{PGcalc_Cint} shows that we can integrate $\fct$ by increasing the subscript of $\GrowthSum{\fct}{r}{\cdot}$ and adjusting the power $r$.
This lemma is a generalization of \cref{PGcalc_int}, which corresponds to the case $j=0$.

\begin{lemma}
\label{PGcalc_Cint}
	Let $\fct \colon \R^d \rightarrow \R^{\dout}$ be $\N_0 \ni k$-times differentiable and $j \in \{0,\dots,k\}$, $r \in [0,\infty)$.
Write $\fct_{\star} \colon \R^{2^jd} \rightarrow \R^{2^j\dout}$ for the $j$-fold augmented derivative of $\fct$.
Then, for any $0 \leq l,m \leq k$ with $j+l+m \leq k$, we have
\begin{equation*}
	\GrowthSum{\nabla^l \fct_{\star}}{r+k-l-m}{m}
	\leq \left( \prod_{i=k-j+1}^k \frac{i(i+1)}{2} \right) \GrowthSum{\fct}{r}{k}.
\end{equation*}%
\end{lemma}

\begin{proof}
	By definition and a shift of summation index, we have
\begin{equation*}
	\GrowthSum{\nabla^l \fct_{\star}}{r+k-l-m}{m}
	= \sup_{x \in \R^d} \sum_{i=l}^{l+m} \frac{\norm{\nabla^i \fct_{\star}(x)}}{(1+\norm{x})^{r+k-i}}
	\leq \GrowthSum{\fct_{\star}}{r+j}{k-j}.
\end{equation*}%
Now, we claim that, for any $\N \ni n$-times differentiable function $\gct \colon \R^{d_1} \rightarrow \R^{d_2}$ and any $s \in [0,\infty)$,
\begin{equation*}
	\GrowthSum{\ED{\gct}}{s+1}{n-1}
	\leq \frac{n(n+1)}{2} \GrowthSum{\gct}{s}{n}.
\end{equation*}%
If this claim holds, then we can apply it repeatedly to find
\begin{equation*}
	\GrowthSum{\fct_{\star}}{r+j}{k-j}
	\leq \left( \prod_{i=k-j}^{k-1} \frac{(i+1)(i+2)}{2} \right) \GrowthSum{\fct}{r}{k}.
\end{equation*}%
It remains to verify the claim.
The case $n=1$ follows from
\begin{equation*}
	\GrowthSum{\ED{\gct}}{s+1}{0}^2
	= \sup_{x,y \in \R^{d_1}} \frac{\norm{\gct(x)}^2 + \norm{\nabla \gct(x)^T y}^2}{(1+\norm{(x,y)})^{2(s+1)}}
	\leq \sup_{x \in \R^{d_1}} \frac{\norm{\gct(x)}^2}{(1+\norm{x})^{2(s+1)}} + \frac{\norm{\nabla \gct(x)}^2}{(1+\norm{x})^{2s}}
	\leq \GrowthSum{\gct}{s}{1}^2.
\end{equation*}%
For a general $n \geq 2$, we use \cref{PGcalc_commute_comDer_stdDer} to estimate
\begin{equation*}
	\GrowthSum{\ED{\gct}}{s+1}{n-1}
	\leq \sup_{x \in \R^{d_1}} \sum_{i=0}^{n-1} \GrowthSum{\nabla^i\ED{\gct}}{s+n-i}{0}
	\leq \sum_{i=0}^{n-1} (i+1) \GrowthSum{\WED{\nabla^i \gct}}{s+n-i}{0}.
\end{equation*}%
Applying the case $n=1$ to the function $\nabla^i \gct$ yields $\GrowthSum{\WED{\nabla^i \gct}}{s+n-i}{0} \leq \GrowthSum{\nabla^i \gct}{s+n-i-1}{1}$.
By \cref{PGcalc_int}, $\GrowthSum{\nabla^i \gct}{s+n-i-1}{1} \leq \GrowthSum{\gct}{s}{n}$, which concludes the proof of the claim.
\end{proof}

When decreasing the subscript of $\GrowthSum{\fct}{r}{\cdot}$, we can get an alternative to \cref{PGcalc_diff} by trading augmented derivatives instead of standard derivatives as shown in the next lemma.

\begin{lemma}
\label{PGcalc_Cdiff}
	Let $\fct \colon \R^d \rightarrow \R^{\dout}$ be $\N_0 \ni k$-times differentiable and $j \in \{0,\dots,k\}$, $r \in [0,\infty)$.
Write $\fct_{\star} \colon \R^{2^jd} \rightarrow \R^{2^j\dout}$ for the $j$-fold augmented derivative of $\fct$.
Then,
\begin{equation*}
	\GrowthSum{\fct}{r}{k}
	\leq \left( \prod_{i=1}^j (1 + 2^r \sqrt{2^id}) \right) \GrowthSum{\fct_{\star}}{r}{k-j}.
\end{equation*}%
\end{lemma}

\begin{proof}
	If we show that, for any $\N \ni n$-times differentiable function $\gct \colon \R^{d_1} \rightarrow \R^{d_2}$,
\begin{equation*}
	\GrowthSum{\gct}{r}{n}
	\leq (1 + 2^r \sqrt{2d_1}) \GrowthSum{\ED{\gct}}{r}{n-1},
\end{equation*}%
then the lemma follows from a repeated application of this inequality.
The factor $2^i$ arises from the fact that the domain of the $i$-fold augmented derivative of $\fct$ is $\R^{2^id}$.
Let us consider the case $n=1$.
Since the Frobenius norm of $\nabla \gct(x)$ can be bounded by the square root of its rank times its operator norm, we can bound
\begin{equation*}
\begin{split}
	\norm{\gct(x)}^2 + \norm{\nabla \gct(x)}^2
	&\leq \norm{\gct(x)}^2 + d_1 \sup_{\norm{y} = 1} \norm{\nabla \gct(x)^T y}^2
	\leq d_1 \sup_{\norm{y} = 1} \norm{\ED{\gct}(x,y)}^2
	\\
	&\leq d_1 \sup_{\norm{y} = 1} \GrowthSum{\ED{\gct}}{r}{0}^2 (1 + \norm{(x,y)})^{2r}
	\leq d_1 \GrowthSum{\ED{\gct}}{r}{0}^2 (2 + \norm{x})^{2r}
\end{split}
\end{equation*}%
and, hence,
\begin{equation*}
	\frac{\norm{\gct(x)}}{(1+\norm{x})^{r+1}} + \frac{\norm{\nabla \gct(x)}}{(1+\norm{x})^r}
	\leq \sqrt{2} \frac{\sqrt{\norm{\gct(x)}^2 + \norm{\nabla \gct(x)}^2}}{(1+\norm{x})^r}
	\leq 2^r \sqrt{2d_1} \GrowthSum{\ED{\gct}}{r}{0}.
\end{equation*}%
Thus, $\GrowthSum{\gct}{r}{1} \leq 2^r \sqrt{2d_1} \GrowthSum{\ED{\gct}}{r}{0}$.
Applying this to $\nabla^n \gct$, we find, for a general $n \geq 2$,
\begin{equation*}
	\GrowthSum{\gct}{r}{n}
	\leq \GrowthSum{\gct}{r+2}{n-2} + \GrowthSum{\nabla^{n-1} \gct}{r}{1}
	\leq \GrowthSum{\gct}{r+2}{n-2} + 2^r \sqrt{2d_1} \GrowthSum{\WED{\nabla^{n-1} \gct}}{r}{0}.
\end{equation*}%
By \cref{PGcalc_commute_comDer_stdDer,PGcalc_Cint}, the second term can be further bound by
\begin{equation*}
	\GrowthSum{\WED{\nabla^{n-1} \gct}}{r}{0}
	\leq \GrowthSum{\nabla^{n-1} \ED{\gct}}{r}{0}
	\leq \GrowthSum{\ED{\gct}}{r}{n-1}.
\end{equation*}%
Finally,
\begin{equation*}
	\GrowthSum{\gct}{r+2}{n-2}
	\leq \GrowthSum{\ED{\gct}}{r+2}{n-2}
	\leq \GrowthSum{\ED{\gct}}{r+1}{n-1}
	\leq \GrowthSum{\ED{\gct}}{r}{n-1},
\end{equation*}%
where the first inequality follows from $\gct(x) = \ED{\gct}(x,0)$, the second one from \cref{PGcalc_Cint}, and the last one from $\GrowthSum{\cdot}{\cdot}{r}$ being non-increasing in $r$.
\end{proof}


\subsection{Polynomial growth of generators}

The next lemma is the first part of \cref{PGcalc_double_generator}.

\begin{lemma}
\label{PGcalc_generator}
	Let $\fct \colon \R^{2d} \rightarrow \R^{\dout}$ be twice differentiable and $r_1,r_2,r_3 \in [0,\infty)$.
Then, for any $t \in [0,T]$,
\begin{equation*}
	\GrowthSum{\GenBO{t}\fct}{r}{0}
	\leq \GrowthSum{\GenBT{t}\fct}{r}{0}
	\leq \max\left\{ \GrowthSum{\mu}{r_1}{0} , \frac{1}{2} \GrowthSum{\sigma}{r_2}{0}^2 \right\} \GrowthSum{\nabla \fct}{r_3}{1}
\end{equation*}%
with $r = \max\{r_1 + 1, 2r_2\} + r_3$.
\end{lemma}

As mentioned after \cref{PGcalc_double_generator}, the same estimate also holds for $\GenO{t}\fct$ and $\GenT{t}\fct$ in place of $\GenBO{t}\fct$ and $\GenBT{t}\fct$ for a function $\fct \colon \R^d \rightarrow \R^{\dout}$ with domain $\R^d$, since we can apply the lemma with $(w,z) \mapsto \fct(w)$.

\begin{proof}
	It is clear that $\GrowthSum{\GenBO{t}\fct}{r}{0} \leq \GrowthSum{\GenBT{t}\fct}{r}{0}$.
Observe that
\begin{equation*}
\begin{split}
	\norm{\GenBT{t}\fct(x,y,z)}
	&\leq \norm{\mu(t,y)} \norm{\nabla \fct(x,z)} + \frac{1}{2} \norm{\sigma(t,y)}^2 \norm{\nabla^2 \fct(x,z)}.
\end{split}
\end{equation*}%
Thus,
\begin{equation*}
\begin{split}
	\frac{\norm{\GenBT{t}\fct(x,y,z)}}{(1+\norm{(x,y,z)})^r} 
	&\leq \frac{\norm{\mu(t,y)}}{(1+\norm{y})^{r_1}} \frac{\norm{\nabla \fct(x,z)}}{(1+\norm{(x,z)})^{r_3+1}} + \frac{1}{2} \frac{\norm{\sigma(t,y)}^2}{(1+\norm{y})^{2r_2}} \frac{\norm{\nabla^2 \fct(x,z)}}{(1+\norm{(x,z)})^{r_3}}
	\\
	&\leq \GrowthSum{\mu}{r_1}{0} \frac{\norm{\nabla \fct(x,z)}}{(1+\norm{(x,z)})^{r_3+1}} + \frac{1}{2} \GrowthSum{\sigma}{r_2}{0}^2 \frac{\norm{\nabla^2 \fct(x,z)}}{(1+\norm{(x,z)})^{r_3}},
\end{split}
\end{equation*}%
from which the desired inequality follows.
\end{proof}

Finally, we fill in the complete proof of \cref{PGcalc_double_generator}.

\begin{proof}[\Pf{PGcalc_double_generator}]
	As mentioned above, the first part of \cref{PGcalc_double_generator} is precisely \cref{PGcalc_generator}.
Now, let us consider the function $\GenBT{t_2}\GenT{t_1}\fct$.
By \cref{PGcalc_generator},
\begin{equation*}
	\GrowthSum{\GenBT{t_2}\GenT{t_1}\fct}{r}{0}
	\leq \max\left\{ \GrowthSum{\mu}{r_1+2}{0} , \frac{1}{2} \GrowthSum{\sigma}{r_2+2}{0}^2 \right\} \GrowthSum{\nabla \GenT{t_1}\fct}{r'}{1},
\end{equation*}%
where
\begin{equation*}
	r' = \max\{r_1+3, 2(r_2+2)\} + r_3 + 2 = r - \max\{r_1 + 3, 2(r_2 + 2)\}.
\end{equation*}%
By \cref{PGcalc_Cint},
\begin{equation*}
	\max\left\{ \GrowthSum{\mu}{r_1+2}{0} , \frac{1}{2} \GrowthSum{\sigma}{r_2+2}{0}^2 \right\}
	\leq \max\left\{ \GrowthSum{\mu}{r_1}{2} , \frac{1}{2} \GrowthSum{\sigma}{r_2}{2}^2 \right\}.
\end{equation*}%
Let $\mu_{\star}$, $\sigma_{\star}$, and $\fct_{\star}$ be the two fold augmented derivatives of $\mu_{t_1}$, $\sigma_{t_1}$, and $\fct$, respectively.
Applying first \cref{PGcalc_Cint}, then \cref{PGcalc_Cdiff}, and then \cref{comDer_generators_commute}, we obtain
\begin{equation*}
	\GrowthSum{\nabla \GenT{t_1}\fct}{r'}{1}
	\leq \GrowthSum{\GenT{t_1}\fct}{r'}{2}
	\leq (1 + 2^{r'} \sqrt{4d}) (1 + 2^{r'} \sqrt{8d}) \GrowthSum{\GenStarT{t_1}\fct_{\star}}{r'}{0},
\end{equation*}%
since $2d$ is the input dimension of $\GenT{t_1}\fct$.
By \cref{PGcalc_generator},
\begin{equation*}
	\GrowthSum{\GenStarT{t_1}\fct_{\star}}{r'}{0}
	\leq \max\left\{ \GrowthSum{\mu_{\star}}{r_1+2}{0} , \frac{1}{2} \GrowthSum{\sigma_{\star}}{r_2+2}{0}^2 \right\} \GrowthSum{\nabla \fct_{\star}}{r_3+2}{1}.
\end{equation*}%
By \cref{PGcalc_Cint},
\begin{equation*}
	\max\left\{ \GrowthSum{\mu_{\star}}{r_1+2}{0} , \frac{1}{2} \GrowthSum{\sigma_{\star}}{r_2+2}{0}^2 \right\}
	\leq \max\left\{ 3 \GrowthSum{\mu}{r_1}{2} , \frac{9}{2} \GrowthSum{\sigma}{r_2}{2}^2 \right\}
\end{equation*}%
and $\GrowthSum{\nabla \fct_{\star}}{r_3+2}{1} \leq 60 \GrowthSum{\fct}{r_3}{4}$.
In conclusion, using that $\max\{a,b/2\} \max\{3a,9b/2\} \leq 9 \max\{a^2,b^2/4\}$,
\begin{equation*}
\begin{split}
	\GrowthSum{\GenBT{t_2}\GenT{t_1}\fct}{r}{0}
	&\leq 540 (1 + 2^{r'} \sqrt{4d}) (1 + 2^{r'} \sqrt{8d}) \max\left\{ \GrowthSum{\mu}{r_1}{2}^2 , \frac{1}{4} \GrowthSum{\sigma}{r_2}{2}^4 \right\} \GrowthSum{\fct}{r_3}{4}
	\\
	&\leq 2^{2r'+12} d \max\left\{ \GrowthSum{\mu}{r_1}{2}^2 , \frac{1}{4} \GrowthSum{\sigma}{r_2}{2}^4 \right\} \GrowthSum{\fct}{r_3}{4}.
\end{split}
\end{equation*}%
The proof for the function $\GenT{t_2}\GenO{t_1}\fct$ is completely analogous with the only exception being the estimate
\begin{equation*}
	\GrowthSum{\nabla \GenO{t_1}\fct}{r'}{1}
	\leq \GrowthSum{\GenO{t_1}\fct}{r'}{2}
	\leq (1 + 2^{r'} \sqrt{2d}) (1 + 2^{r'} \sqrt{4d}) \GrowthSum{\GenStarO{t_1}\fct_{\star}}{r'}{0},
\end{equation*}%
since the input dimension of $\GenO{t_1}\fct$ is $d$ whereas before we had a factor of $2d$.
\end{proof}

\end{appendix}


\vskip 5mm\noindent{\large\scshape Acknowledgments}\vskip 2mm
\noindent
We are grateful to Philipp Zimmermann and Robert Crowell for fruitful discussions and helpful comments.

\phantomsection%
\addcontentsline{toc}{section}{\protect Bibliography}%
\bibliographystyle{acm}%
\bibliography{bibfile_FR_2023_june_9}

\end{document}